\documentclass[hidelinks,onefignum,onetabnum]{siamart250211}

\usepackage{amsmath,amsfonts,amssymb}
\usepackage{xcolor,bm,url}
\usepackage{booktabs}
\usepackage{tikz}
\usepackage{hyperref}
\usepackage{etoolbox}
\usepackage{xspace}

\newsiamremark{remark}{Remark}

\DeclareMathOperator*{\argmax}{arg\,max}

\def\nlp{NLP($p$)}
\def\lagrangian{\mathcal{L}}

\def\activeset{\mathcal{I}}
\def\multiplierset{\mathcal{M}}
\def\solutionset{\Sigma}
\def\extremeset{\mathcal{E}}
\def\criticalcone{\mathcal{K}}
\def\nullspace{N}
\def\chull{\text{conv}}
\def\valuefunction{\varphi}

\headers{Sensitivity analysis for parametric NLP}{F. Pacaud}

\title{Sensitivity analysis for parametric nonlinear programming: A tutorial}
\author{François Pacaud \thanks{Centre Automatique et Systèmes (CAS), Mines Paris-PSL.}}
\date{\today}

\begin{document}
\maketitle
\begin{abstract}
  This tutorial provides an overview of the current state-of-the-art in
  the sensitivity analysis for nonlinear programming.
  Building upon the fundamental work of Fiacco, it derives the sensitivity
  of primal-dual solutions for regular nonlinear programs and explores the extent
  to which Fiacco's framework can be extended to degenerate nonlinear programs
  with non-unique dual solutions.
  The survey ends with a discussion on how to adapt the sensitivity analysis for conic programs
  and approximate solutions obtained from interior-point algorithms.
\end{abstract}

\section{Introduction}

Sensitivity analysis has been a long standing research thread in
nonlinear programming. It aims at evaluating the
sensitivity of a solution $x \in \mathbb{R}^n$ of a nonlinear problem with respect to parameter changes.
For a parameter $p\in\mathbb{R}^\ell$, the solution of a generic optimization problem
is given as solution of a (possibly nonsmooth) variational problem:
\begin{equation}
  \label{eq:system}
  F(x, p) = 0 \; ,
\end{equation}
with $F: \mathbb{R}^n \times \mathbb{R}^\ell \to \mathbb{R}^n$ a function encoding the stationary conditions
of the optimization problem (e.g. the Karush-Kuhn-Tucker stationary conditions,
or the homogeneous self-dual embedding associated to a conic program).
The solution of \eqref{eq:system} depends on $p$ and is noted $x(p)$.
If the function $F$ is \emph{regular} enough, the Implicit Function Theorem
gives an elegant way to evaluate the Jacobian $J_p x(p) \in \mathbb{R}^{n \times \ell}$ of the solution $x(p)$ of \eqref{eq:system} as
the solution of the linear system:
\begin{equation}
  \label{eq:iftderiv}
  J_p x(p) = - \big(J_x F(x, p)\big)^{-1} J_p F(x, p)  \; .
\end{equation}
Most applications require only the forward and adjoint sensitivities,
alleviating the need to compute the full Jacobian $J_p x(p)$.
For vectors $u \in \mathbb{R}^\ell$ and
$v \in \mathbb{R}^n$, the forward and adjoint sensitivities are defined respectively as
\begin{equation}
  \begin{aligned}
    & J_p x(p)\phantom{^\top} \, u = - \big(J_x F(x, p)\big)^{-1} \, \bar{u} && \text{with } \bar{u} := J_p F(x, p) \, u  \; , \\
    & J_p x(p)^\top \, v = - J_p F(x, p)^\top \, \bar{v} && \text{with } \bar{v} := \big(J_x F(x, p)\big)^{-\top} \, v   \; .
  \end{aligned}
\end{equation}

Following \cite{blondel2022efficient},
we call \emph{implicit differentiation} the class of differentiation methods based on the relation \eqref{eq:iftderiv}.
Implicit differentiation requires solving the nonlinear system~\eqref{eq:system}
exactly and  then relies on the Implicit Function Theorem: it is accurate,
but often limited in scope. A practical alternative uses the fact that
the system~\eqref{eq:system} is usually solved \emph{iteratively}, e.g. using an optimization
solver. If the iterative method is differentiable (in the sense each of its operation can
be differentiated by automatic differentiation), the sensivitity $J_p x(p)$ can
be computed \emph{iteratively} using forward accumulation (also
known as \emph{piggy back} method~\cite{griewank2008evaluating}). This second mode does not rely on the Implicit Function
Theorem, as it propagates the derivatives directly through the algorithm. The main advantage is
the derivatives are directly computed at the solution, meaning we do not have to solve
an additional linear system as in \eqref{eq:iftderiv}. However, the method
can be sensitive to numerical error and the derivatives $J_p x(p)$ converge
more slowly than the primal solution $x(p)$ \cite{gilbert1992automatic,griewank2002reduced}.

Hence, implicit differentiation remains the favorite method to compute the sensitivities.
Unfortunately, the hypothesis of the Implicit Function Theorem are often too stringent
to hold for a generic nonlinear program, limiting the applicability of \eqref{eq:iftderiv}
unless practical workarounds are applied.
For nonlinear programs, the problem's regularity
is related to its non-degeneracy: if we perturb slightly the parameters,
we have to ensure that the primal solution exists and remains unique.
Depending on the structure of the problem, the perturbed solution may
have the following demonstrable properties~\cite{fiacco1993degeneracy}, ordered by difficulty:
\begin{itemize}
  \item Existence of the primal solution.
  \item Existence and continuity of the primal solution.
  \item (Upper) Lipschitz continuity of the primal solution, non-uniqueness of the primal-dual solution.
  \item (Upper) Lipschitz continuity of the primal-dual solution, non-uniqueness of the primal-dual solution.
  \item Directional differentiability and uniqueness of the primal solution, non-uni\-queness of the dual solution.
\end{itemize}
In this tutorial, we aim at summarizing the conditions under which the sensitivity
analysis is well-posed.

\subsection{A brief history of sensitivity analysis for nonlinear programming}
We start by a brief review of literature to remind the
main developments in the theory of sensitivity analysis for nonlinear programming.

\paragraph{The 1970s: defining regularity and differentiability}
The regularity of the KKT solution has been studied extensively in the 1970s
by Robinson~\cite{robinson1980strongly,robinson1982generalized}, by interpreting
the nonlinear program as a generalized nonlinear equation.
It is now established that the regularity of a nonlinear program usually follows from
the constraint qualification conditions satisfied at the KKT solution.
Robinson proved in \cite{robinson1982generalized} that the primal-dual solution
is \emph{strongly regular} (hence, locally unique and Lipschitz continuous)
if LICQ and SSOSC holds at the original solution.
In addition, the (Fréchet) differentiability of the primal-dual solution is established
if SCS hold, as established by Fiacco in his seminal result~\cite{fiacco1976sensitivity}.

\paragraph{The 1980s and 1990s: tackling degeneracy}
The previous results have been extended in the 1980s, with a primary focus
on dropping the strong SCS condition required in~\cite{fiacco1976sensitivity}.
Kojima~\cite{kojima1980strongly} proved that MFCQ and GSSOSC are sufficient
to guarantee that the primal-dual solution is strongly stable.
Moreover, the primal-dual solution becomes non-differentiable as soon as we drop the
SCS condition, but it remains \emph{directionally differentiable} under certain
assumptions.
Jittorntrum~\cite{jittorntrum1984solution} showed that if LICQ and SSOSC hold,
then the directional derivative of the primal-dual solution exists and is
given as the solution of a quadratic program (QP).

If we drop the LICQ condition, the dual solution becomes non-unique and
the directional differentiability of the dual solution is lost.
Shapiro~\cite{shapiro1985second} showed that the primal solution
remains directionally differentiable if we replace LICQ by SMCFQ and SSOSC. This
result was extended by Kyparisis in \cite{kyparisis1990sensitivity}, who proved
that we retain directionally differentiability for the primal solution if we assume MFCQ, CRCQ and
GSSOSC, that is, without requiring unicity of the dual solution. This thread of
research culminated in 1995 with the publication of the seminal result of Ralph and
Dempe~\cite{ralph1995directional}, giving a practical method to evaluate
the directional derivative of the primal solution under MFCQ.

\paragraph{The 2000s and 2010s: extending the scope of applications}
The theory of sensitivity analysis for nonlinear program has witnessed
a significant consolidation after the 1990s, with the publication
of three consequent books dedicated to the sensitivity analysis of optimization programs~\cite{bonnans2013perturbation,facchinei2003finite,dontchev2009implicit}.
In particular, the application of sensitivity analysis has been extended to the more
general setting offered by variational inequalities, which covers both nonlinear programs
and complementarity problems.
In parallel, numerous applications have flourished out of the sensitivity
analysis theory, one of the most significant being
bilevel programming~\cite{dempe2002foundations}.

\subsection{Motivations}
Sensitivity analysis has proven to be a fertile ground for
the analysis of stochastic optimization problems~\cite{shapiro1990differential,shapiro1991asymptotic},
parametric decomposition~\cite{gauvin1978method,demiguel2008decomposition},
parameter estimation problems~\cite{lopez2012moving}
and last but not least, model predictive control~\cite{zavala2009advanced,jaschke2014fast,zavala2010real}.
Recently, the field has witnessed a renewed interest with the emergence
of the \emph{differentiable programming} paradigm in the machine learning community,
where people look at embedding a (convex) optimization solver inside a neural network
\cite{amos2017optnet,agrawal2019differentiable} (whose training requires the evaluation of the sensitivity in the backward pass).

\subsection{Outline}
Section~\ref{sec:nlp} is devoted to the theory of nonlinear
programming, with an emphasis on the constraint qualifications used to define
the regularity of the feasible set.
Section~\ref{sec:solutionsensitivity} summarizes the main results
about the sensitivity of the primal-dual solution in nonlinear programming.
It extends Fiacco's seminal work \cite{fiacco1983introduction}
and puts a special consideration on degenerate nonlinear programs, following
the lead of Ralph and Dempe~\cite{ralph1995directional}.
Section~\ref{sec:objectivesensitivity} focuses more particularly on
the sensitivity of the objective function. Last, Section~\ref{sec:numerics}
details practical methods and specific software implementations to compute the sensitivity
of any optimization program, including an extended discussion on the
sensitivity of conic programs.

\section{Nonlinear programming}
\label{sec:nlp}
This section recalls classical results in the theory of nonlinear programming,
with a long discussion on constraint qualifications \S\ref{sec:nlp:cq}.
The scope of this tutorial is limited to finite-dimensional optimization problems.
We refer to \cite{bonnans2013perturbation} for a broader description covering the
infinite dimensional cases.

\subsection{Notations}
For a parameter $p \in \mathbb{R}^{\ell}$, we are
interested in solving the parametric nonlinear optimization problem:
\begin{equation}
  \label{eq:problem}
  \text{NLP}(p) := ~
  \min_{x \in \mathbb{R}^{n}} \; f(x, p) \quad \text{subject to} \quad
  \left\{
  \begin{aligned}
    g(x, p) = 0 \;,\\
    h(x, p) \leq 0 \;,
  \end{aligned}
  \right.
\end{equation}
where $f:\mathbb{R}^{n} \times \mathbb{R}^{\ell} \to \mathbb{R}$ an
objective and
$g:\mathbb{R}^{n} \times \mathbb{R}^{\ell} \to \mathbb{R}^{m_e}$ and
$h:\mathbb{R}^{n} \times \mathbb{R}^{\ell} \to \mathbb{R}^{m_i}$ two functions
encoding respectively the equality and inequality constraints.

All functions are twice-differentiable and depend jointly on the optimization variable $x \in \mathbb{R}^{n}$
and a parameter $p \in \mathbb{R}^{\ell}$.
The \emph{feasible set} of the optimization problem~\eqref{eq:problem} is encoded by the set-valued function
\begin{equation}
  X(p) = \{x \in \mathbb{R}^{n} \; | \; g(x, p) = 0 \, , \, h(x, p) \leq 0  \} \, .
\end{equation}

\paragraph{Value function}
We define the value function $\valuefunction:\mathbb{R}^{\ell} \to
\mathbb{R} \cup \{+\infty\}$ as the solution of \nlp\
for a given parameter $p \in \mathbb{R}^{\ell}$:
\begin{equation}
  \label{eq:valuefunction}
  \valuefunction(p) = \min_{x \in \mathbb{R}^{n}} \; f(x, p) \quad \text{subject to} \quad
  x \in X(p) \;.
\end{equation}
By convention $\valuefunction(p) = +\infty$ if the problem is infeasible ($X(p) = \varnothing$).
The \emph{solution set} associated to \eqref{eq:valuefunction} is noted
\begin{equation}
  \label{eq:solutionset}
  \solutionset(p) := \{x \in X(p)  \; : \; f(x, p) = \valuefunction(p) \} \;.
\end{equation}

\paragraph{Lagrangian}
We introduce the Lagrangian associated to \eqref{eq:problem} as
\begin{equation}
  \label{eq:lagrangian}
  \lagrangian(x, y, z, p) = f(x, p) + y^\top g(x, p) + z^\top h(x, p) \;,
\end{equation}
where $y \in \mathbb{R}^{m_e}$ (resp. $z \in \mathbb{R}^{m_i}$) the multiplier associated to the equality
(resp. inequality) constraints.
The Karush-Kuhn-Tucker (KKT) conditions of~\eqref{eq:problem}
state that under a certain constraint qualification, the primal vector $x$
is a local optimal solution if there exists a dual multiplier $(y, z)$ such that
\begin{equation}
  \label{eq:kkt}
  \begin{aligned}
    \nabla_x f(x, p) + \sum_{i=1}^{m_e} y_i \nabla_x g_i(x, p) + \sum_{i=1}^{m_i} z_i \nabla_x h_i(x, p) = 0 \;, \\
    g(x, p) = 0 \, , \quad h(x, p) \leq 0 \,, \quad z \geq 0 \, ,\quad z^\top h(x, p) = 0 \;.
  \end{aligned}
\end{equation}
The KKT conditions become necessary and sufficient in the convex case.
We define a primal-dual stationary solution (or KKT stationary solution) as the vector
\begin{equation}
  w(p) = \big(x(p), y(p), z(p) \big)  \in \mathbb{R}^{n} \times \mathbb{R}^{m_e} \times \mathbb{R}^{m_i} \; .
\end{equation}
To alleviate the notations, we often denote by $w^\star$ the vector $w(p)$.


\paragraph{Multiplier set}
The set of multipliers $(y, z) \in \mathbb{R}^{m_e+m_i}$ satisfying the KKT conditions \eqref{eq:kkt}
is denoted by $\multiplierset_p(x)$. We notice that the set $\multiplierset_p(x)$ is polyhedral (hence closed
and convex).
We note by $\extremeset_p(x)$ the set of extreme points of the polyhedron $\multiplierset_p(x)$.
We say that the problem is \emph{degenerate} if $\multiplierset_p(x)$ has more than one element.

\paragraph{Active set}
For a given $x \in \mathbb{R}^{n}$, we define the \emph{active set} at $x$ as
\begin{equation}
  \activeset_p(x) = \{ i \in [m_i] \;|\; h_i(x, p) =0 \}
  \; .
\end{equation}
For given multipliers $(y, z) \in \multiplierset_p(x)$, we define the \emph{strongly} and
\emph{weakly active constraints} as
\begin{equation}
  \activeset^+_p(x; z) = \{ i \in \activeset_p(x) \; | \; z_i > 0 \} \;, \quad
  \activeset^0_p(x; z) = \{ i \in \activeset_p(x) \; | \; z_i = 0 \} \; .
\end{equation}
The strongly and weakly active constraints form a partition of the active set:
$\activeset_p(x) = \activeset_p^+(x; z) \cup \activeset_p^0(x; z)$
and $\activeset_p^+(x; z) \cap \activeset_p^0(x; z) = \varnothing$.

Degeneracy arises from the undetermined weakly active constraints $\activeset_p^0(x, z)$:
a small perturbation of the parameter $p$ is likely to render them inactive.
We say that the problem satisfies \emph{strict complementarity slackness} (SCS) if
there is no weakly active constraint: $\activeset_p^0(x; z) = \varnothing$.

\paragraph{Critical cone}
The critical cone is a key notion to characterize the second-order
optimality conditions of a nonlinear program.

\begin{definition}[Critical cone]
  For a feasible point $x \in X(p)$,
  the critical cone associated to the feasible set $X(p)$ at $x$ is defined as
  \begin{equation}
    \label{eq:criticalcone}
    \begin{aligned}
      \mathcal{C}_p(x) = \Big\{ d \in \mathbb{R}^{n} \; : \;
          & \nabla_{x} f(x, p)^\top d \leq 0 \;, \\
          & \nabla_{x} g_i(x, p)^\top d = 0 \;,\forall i \in [m_e] \;, \\
          & \nabla_{x} h_j(x, p)^\top d \leq 0 \; , \forall j \in \activeset_p(x) \Big\} \; .
    \end{aligned}
  \end{equation}
\end{definition}
We note that at a KKT solution $(x^\star, y^\star, z^\star)$, the critical cone simplifies as
\begin{equation}
    \begin{aligned}
      \mathcal{C}_p(x^\star) = \Big\{ d \in \mathbb{R}^{n} \; : \;
          & \nabla_{x} g_i(x^\star, p)^\top d = 0 \;,\forall i \in [m_e] \;, \\
          & \nabla_{x} h_j(x^\star, p)^\top d = 0 \; , \forall j \in \activeset_p^+(x^\star; z^\star) \;, \\
          & \nabla_{x} h_j(x^\star, p)^\top d \leq 0 \; , \forall j \in \activeset_p^0(x^\star; z^\star) \Big\} \; .
    \end{aligned}
\end{equation}
Further, if we assume that strict complementarity slackness holds, the critical cone is exactly the null space
of the active constraints' Jacobian: $\mathcal{C}_p(x^\star) = \nullspace(J_{act}(x^\star))$:
\begin{equation}
  \begin{aligned}
    \mathcal{C}_p(x^\star) = \{ d \in \mathbb{R}^{n} \; : \; &
      \nabla_{x} g_i(x^\star, p)^\top d = 0 \;, \forall i \in [m_e] \,, \\
      & \nabla_{x} h_j(x^\star, p)^\top d = 0 \; , \forall j \in \activeset_p(x^\star) \} \; .
  \end{aligned}
\end{equation}
In parametric optimization, it is common to introduce the critical set
associated to the critical cone~\eqref{eq:criticalcone} at the solution~\cite{facchinei2003finite,ralph1995directional}.
\begin{definition}[Critical set]
  Let $(x^\star, y^\star, z^\star)$ be a KKT stationary solution.
  The \emph{critical set}
  at the solution is the polyhedral cone defined as
  \begin{equation}
    \begin{aligned}
      \criticalcone_p(x^\star, z^\star) = \{ d \in \mathbb{R}^{n + \ell} \; : \;
          & \nabla_{(x,p)} g_i(x^\star, p)^\top d = 0 \;, \forall i \in [m_e] \,,\\
          & \nabla_{(x,p)} h_j(x^\star, p)^\top d = 0 \; , \forall j \in \activeset_p^+(x^\star; z^\star) \,, \\
          & \nabla_{(x,p)} h_j(x^\star, p)^\top d \leq 0 \; , \forall j \in \activeset_p^0(x^\star; z^\star) \} \; .
    \end{aligned}
  \end{equation}
  Elements of $\criticalcone_p(x^\star, z^\star)$ are called \emph{critical directions}.
  The \emph{directional critical set} in the direction $h \in \mathbb{R}^{\ell}$
  is defined as
  \begin{equation}
    \criticalcone_p(x^\star, z^\star; h) = \{ d \in \mathbb{R}^{n} \; : \; (d, h) \in \criticalcone_p(x^\star, z^\star) \} \;.
  \end{equation}
\end{definition}

\subsection{Constraint qualifications}
\label{sec:nlp:cq}
Constraint qualifications ensure that we can capture \emph{locally} the geometry
of the feasible set $X(p)$ by linearizing it near a feasible point $x \in X(p)$.
The regularity of a nonlinear program is defined by the particular constraint qualifications
that hold at a stationary point. We recall in this section the main constraint
qualifications we can encounter in practice, and their respective impacts on the problem's topology.

\begin{definition}[Linear-independence constraint qualification]
  LICQ holds if the gradients of active constraints
  \begin{equation}
  \{\nabla_x g_i(x, p) \; : \; i \in [m_e]\}~\cup~
  \{\nabla_x h_i(x, p) \; : \; i \in \activeset_p(x) \} \; ,
  \end{equation}
  are linearly independent.
\end{definition}

LICQ implies the following second-order necessary condition.

\begin{theorem}
  \label{thm:soclicq}
  Suppose $(x^\star, y^\star, z^\star)$ is a primal-dual solution
  of \eqref{eq:problem} satisfying LICQ. Then,
  \begin{equation}
    d^\top \nabla^2_{x x} L(x, y, z, p) d \geq 0 \quad \forall
    d \in \mathcal{C}_p(x) \; ,
  \end{equation}
  where $\mathcal{C}_p(x)$ is the critical cone of the feasible set $X(p)$ at $x$.
\end{theorem}

\begin{definition}[Mangasarian-Fromovitz constraint qualification]
  MFCQ holds if the gradients of the equality constraints
  \begin{equation}
    \{\nabla_x g_i(x, p) \; : \; i \in [m_e]\} \; ,
  \end{equation}
  are linearly independent
  and if there exists a vector $d \in \mathbb{R}^{n_x}$ such that
  \begin{equation}
      \nabla_x g_i(x, p)^\top d = 0 \;,  \forall i \in [m_e] \, , \quad
    \nabla_x h_i(x, p)^\top d < 0 \;, \forall i \in \activeset_p(x) \, .
  \end{equation}
\end{definition}
It is well known that MFCQ generalizes the Slater's condition of convex
optimization to nonlinear programming. MFCQ also implies the boundedness
of the optimal multipliers set.

\begin{theorem}[\cite{gauvin1977necessary}]
  The set $\multiplierset_p(x)$ is bounded if and only if MFCQ holds.
\end{theorem}

\begin{proposition}[Dual MFCQ]
  Using the theorem of alternatives, MFCQ is
  equivalent to $0$ being the unique solution of the linear system
  \begin{equation}
    \begin{aligned}
    & \sum_{i=1}^{m_e} y_i \nabla_x g_i(x, p) +
    \sum_{i\in \activeset_p(x)} z_i \nabla_x h_i(x, p) = 0 \; , \\
    & z_i \geq 0  \quad \forall i \in \activeset_p(x) \; .
    \end{aligned}
  \end{equation}
\end{proposition}

\begin{definition}[Strict Mangasarian-Fromovitz constraint qualification] \\
  SMFCQ holds if
  \begin{equation}
  \{\nabla_x g_i(x, p) \; : \; i \in [m_e]\}~\cup~
  \{\nabla_x h_i(x, p) \; : \; i \in \activeset_p^+(x; z) \} \;,
  \end{equation}
  are linearly independent and if there exists a vector $d \in \mathbb{R}^{n_x}$ such that
  \begin{equation}
    \begin{aligned}
    & \nabla_x g_i(x, p)^\top d = 0 \;, \forall i \in [m_e] \, , \\
    & \nabla_x h_i(x, p)^\top d = 0 \;, \forall i \in \activeset_p^+(x; z) \, , \\
    & \nabla_x h_i(x, p)^\top d < 0 \;, \forall i \in \activeset_p^0(x; z) \, .
    \end{aligned}
  \end{equation}
\end{definition}
A primal-dual solution $(x^\star, y^\star, z^\star)$ satisfying SMCFQ
implies a second-order necessary condition analogous to Theorem \ref{thm:soclicq}
\cite[Theorem 2.1]{kyparisis1985uniqueness}.
SMFCQ is the weakest condition under which the multiplier set $\multiplierset_p(x)$ reduces to a singleton.
\begin{proposition}[Proposition 1.1, \cite{kyparisis1985uniqueness}]
  \label{prop:kyparisis1985}
  Let $w = (x, y, z)$ be a primal-dual solution. Then SMFCQ
  holds at $x$ if and only if the multiplier vector $(y, z)$ is unique.
\end{proposition}

\begin{definition}[Constant-rank constraint qualification (CRCQ)]
  CRCQ holds if there exists a neighborhood $W$ of $x^\star$ such that for
  any subsets $I$ of $\activeset_p(x^\star)$ and $J$ of $[m_e]$ the family
  of gradient vectors
  \begin{equation}
  \{\nabla_x g_j(x, p) \; : \; j \in J\}~\cup~ \{\nabla_x h_i(x, p) \; : \; i \in I \}
  \; ,
  \end{equation}
  has the same rank for all vectors $x \in W$.
\end{definition}

We now focus on second-order constraint qualifications, starting
with the second-order sufficiency condition.

\begin{definition}[Second-order sufficiency condition (SOSC)]
  SOSC holds at $x$ if there exists $(y, z) \in \multiplierset_p(x)$ such that
  \begin{equation}
    d^\top \nabla^2_{x x} L(x, y, z, p) d > 0 \quad \forall
    d \in \mathcal{C}_p(x) \;, \;d \neq 0 \;.
  \end{equation}
\end{definition}

\begin{proposition}
  Suppose SOSC holds at a solution $x^\star$.
  Then, there exists $\sigma > 0$ and a neighborhood $V$ of
  $x^\star$ such that for all $x \in X(p) \cap V$
  (with $x \neq x^\star$),
  \begin{equation}
    f(x) > f(x^\star) + \frac{\sigma}{2} \|x - x^\star \|^2 \; .
  \end{equation}
\end{proposition}

We define the set
\begin{equation}
  \begin{aligned}
  \mathcal{D}_p(x; z) = \{ d \in \mathbb{R}^{n_x} \;:\;
  & \nabla_x g_i(x, p) d = 0 \quad \forall i \in [m_e] \;, \\
  & \nabla_x h_j(x, p) d = 0  \quad \forall j \in \activeset_p^+(x; z) \} \;.
  \end{aligned}
\end{equation}
Note that at a primal-dual solution $(x, y, z)$, $\mathcal{C}_p(x; z) \subset \mathcal{D}_p(x; z)$
with both sets being equal if SCS holds.

\begin{definition}[Strong second-order sufficiency condition (SSOSC)]
  SSOSC holds at $x$ if there exists $(y, z) \in \multiplierset_p(x)$ such that
  \begin{equation}
    \label{eq:cq:ssosc}
    d^\top \nabla^2_{x x} L(x, y, z, p) d > 0  \; \quad \forall
    d \in \mathcal{D}_p(x; z) \;, d \neq 0 \;.
  \end{equation}
  If the condition \eqref{eq:cq:ssosc} hold for all multipliers in $\multiplierset_p(x)$,
  then we say that the generalized second-order sufficiency condition (GSSOSC) hold.
\end{definition}
If SCS hold, we note that SSOSC becomes equivalent to SOSC.
\section{Sensitivity of the primal-dual solution}
\label{sec:solutionsensitivity}
Let $p^\star \in \mathbb{R}^{\ell}$ be a parameter, and
$w^\star = (x^\star, y^\star, z^\star)$ a KKT stationary solution of \nlp.
We aim at studying the behavior of the primal-dual
solution $w(p) := (x(p), y(p), z(p))$ near $p^\star$.
In particular, we are interested in characterizing the local continuity of the solution $w(p)$
(i.e., \emph{is $w(p)$ Lipschitz-continuous?}) as well as its local differentiability
(i.e., \emph{under which conditions is $w(p)$ differentiable?}).

We start by recalling Fiacco's theorem in \S\ref{sec:solutionsensitivity:fiacco}.
This result is often restrictive, as it does not allow for active set changes
in the solution. For that reason, the sensitivity are often derived
using Robinson's generalized equations, as described in \S\ref{sec:solutionsensitivity:robinson}.
The generalized equation framework is a convenient way to encode the nonsmoothness
in the problem's solution. In \S\ref{sec:solutionsensitivity:lexicographic}, we show
that it can also be interpreted using the notion of lexicographic derivative.
The last section \S\ref{sec:solutionsensitivity:degenerate} is devoted
to the differentiation of degenerate nonlinear program, where the dual solution
is not unique. We show that in that case the notion of $PC^1$ functions
is the most adapted to characterize the nature of the solution.

\subsection{Fiacco and the implicit function theorem}
\label{sec:solutionsensitivity:fiacco}
If the problem is regular, the local continuity is usually given by
the Implicit Function Theorem.
\begin{theorem}[Implicit function theorem]
  \label{thm:ift}
  Let $F: \mathbb{R}^{n} \times \mathbb{R}^{\ell} \to \mathbb{R}^{n}$ be a $C^1$ function
  and $(x^\star, p^\star) \in \mathbb{R}^{n} \times \mathbb{R}^{\ell}$ such that $F(x^\star, p^\star) = 0$.
  If the Jacobian $J_x F(x^\star, p^\star)$ is \emph{invertible}, then there exists
  a neighborhood $U \subset \mathbb{R}^{\ell}$ of $p^\star$ and a differentiable
  function $x(p)$ such that $F(x(p), p)= 0$ for all $p \in U$,
  and $x(p)$ is the unique solution in a neighborhood of $x^\star$.
  In addition, $x(\cdot)$ is $C^1$ differentiable on $U$, with
\begin{equation}
  \label{eq:ift}
  J_p x(p) = -\big(J_x F(x, p)\big)^{-1} J_p F(x, p)
  \; .
\end{equation}
\end{theorem}
In \cite{fiacco1976sensitivity}, Fiacco investigated under which conditions the KKT equations falls
under the conditions of the implicit function theorem.
To do so, we should ensure (i) the primal solution mapping $x(p)$ is
single-valued (i.e. the problem \eqref{eq:problem} has an unique solution)
(ii) the dual solution $(y(p), z(p))$ is unique and (iii) the active set is locally stable.
Fiacco showed that these three conditions are satisfied if respectively
(i) SOSC (ii) LICQ and (iii) SCS hold.

In that particular case, the KKT equations~\eqref{eq:kkt} rewrite as a smooth system
of nonlinear equations depending continuously on $w = (x, y, z)$.
For $i \in \activeset_p(x)$ we have
$z_i > 0$ and $h_i(x, p) = 0$;
alternatively for $i \in [m_i] \setminus \activeset_p(x)$ we have
$z_i = 0$ and $h_i(x, p) < 0$.
By differentiating the complementarity condition $z_i h_i(x, p) = 0$,
we get for all $i = 1,\cdots, m_i$,
\begin{equation}
  \dfrac{\partial z_i}{\partial p} \cdot h_i(x, p) +
  z_i \cdot \left( \dfrac{\partial h_i}{\partial x}\dfrac{\partial x}{\partial p} +
  \dfrac{\partial h_i}{\partial p} \right) = 0 \; .
\end{equation}
Hence, if $i \in [m_i] \setminus \activeset_p(x)$ we get by elimination $\frac{\partial z_i}{\partial p} = 0$:
the multiplier $z_i = 0$ remains constant after a small perturbation, meaning the
constraints remain inactive.
Alternatively, if $i \in \activeset_p(x)$, we deduce that
$\frac{\partial h_i}{\partial x}\frac{\partial x}{\partial p} + \frac{\partial h_i}{\partial p} = 0$
: the constraint remains locally active.
As a consequence, we can discard the inactive constraints information
from the KKT conditions and interpret the remaining active constraints as equality
constraints. We note $z_\mathcal{I} = \{z_i \}_{i \in \activeset_p}$ the vector
storing the active multipliers and $h_\mathcal{I}(x, p) = \{h_i(x, p) \}_{i \in \activeset_p}$
the vector storing the active inequality constraints.

Locally, the solution $w := (x, y, z_\mathcal{I})$ satisfies the smooth nonlinear system:
\begin{equation}
  \label{eq:kktlagr}
  F(w) = \begin{pmatrix}
    \nabla_x \mathcal{L}(x, y, z_\mathcal{I}) \\
    g(x, p) \\
    h_\mathcal{I}(x, p)
  \end{pmatrix}
  = 0 \; .
\end{equation}
We note the Jacobian of $F(\cdot)$ :
\begin{equation}
  M := J_w F =
\begin{bmatrix}
  \nabla^2_{x x} \lagrangian & J_{x} g^\top & J_{x} h_\mathcal{I}^\top \\
  J_x g & 0 & 0 \\
  J_x h_\mathcal{I} & 0 & 0
\end{bmatrix}
\;, \quad
N := J_p F =
\begin{bmatrix}
  \nabla^2_{x p} \lagrangian \\
  J_p g  \\
  J_p h_\mathcal{I}
\end{bmatrix}
\; .
\end{equation}

The sensitivity of the solution is given by the following theorem,
which applies the Implicit Function Theorem~\ref{thm:ift} to the
system \eqref{eq:kktlagr}.

\begin{theorem}[Theorem 3.2.2 \cite{fiacco1976sensitivity}]
  \label{thm:fiacco}
  Suppose that at a KKT stationary solution $x^\star$ of NLP$(p^\star)$, SOSC, LICQ and SCS hold. Then,
  \begin{enumerate}
    \item  the primal-dual solution $w^\star = (x^\star, y^\star, z^\star_\mathcal{I})$ is unique;
    \item  there exists a unique continuously differentiable function $w(\cdot)$
      defined in a neighborhood $U$ of $p^\star$ such that $w(p^\star) =(x^\star, y^\star, z^\star_\mathcal{I})$, with
  \begin{equation}
    \label{eq:fiaccoderiv}
    J_p w(p) = - M^{-1} N \; .
  \end{equation}
\item SCS and LICQ hold locally near $p^\star$.
  \end{enumerate}
\end{theorem}
In detail, Equation~\eqref{eq:fiaccoderiv} gives the sensitivity
of the primal-dual solution as
\begin{equation}
  \label{eq:fiaccoderiv2}
  \begin{bmatrix}
  J_p x  \\
  J_p y  \\
  J_p z_\mathcal{I}
  \end{bmatrix}
 = - \begin{bmatrix}
  \nabla^2_{x x} \lagrangian & J_{x} g^\top & J_{x} h_\mathcal{I}^\top \\
  J_x g & 0 & 0 \\
  J_x h_\mathcal{I} & 0 & 0
\end{bmatrix}^{-1}
\begin{bmatrix}
  \nabla^2_{x p} \lagrangian \\
  J_p g  \\
  J_p h_\mathcal{I}
\end{bmatrix}
  \; .
\end{equation}
The operation translates as the solution of a (sparse) linear system with $\ell$ right-hand-sides.
We note that the linear system \eqref{eq:fiaccoderiv2} has a saddle point structure:
the SOSC and LICQ conditions implies that it is invertible~\cite[Theorem 3.4]{benzi2005numerical}.
If the original problem \nlp\ is solved with a Newton algorithm, the factorization of the Jacobian matrix
$J_w F$ is already available, leaving only to compute $\ell$ backsolves.

\begin{remark}[Linear programming]
  If the problem~\eqref{eq:problem} is a linear program (LP), we have
  $\nabla_{xx}^2 \mathcal{L} = 0$. The Goldman-Tucker theorem
  states there always exists a LP solution satisfying strict-complementarity. If
  the LP is not degenerate, LICQ holds and the solution is a vertex
  of the polyhedral feasibility set $X(p)$. This automatically implies SOCP (the critical cone
  reduces to $\mathcal{C}_p(x^\star) = \{0 \}$) and there are exactly
  $n$ binding constraints: $m_e + | \activeset_p(x^\star) | = n$.
  Hence Theorem~\ref{thm:fiacco} applies to the linear programming case.
  We denote the active Jacobian at the solution $B = \begin{bmatrix}
    J_x g^\top & J_x h_\mathcal{I}^\top
  \end{bmatrix}^\top$.
  The matrix $B$ has dimension $n \times n$ and is non-singular (LICQ implies
  it is full row-rank). We note that $B$ is exactly the basis matrix used in the simplex algorithm.
  In that case, the inverse of the matrix $M = \begin{bmatrix} 0 & B^\top \\ B & 0 \end{bmatrix}$
  in \eqref{eq:fiaccoderiv} satisfies:
  \begin{equation}
    M^{-1} = \begin{bmatrix}
      0 & B^{-1} \\ B^{-\top} & 0
    \end{bmatrix} \; .
  \end{equation}
  Hence, the computation of the sensitivities simplifies considerably
  if the problem is regular linear program. If the linear program is solved
  using the simplex algorithm, the solver usually returns a LU factorization
  for the basis matrix $B$.
\end{remark}

\begin{remark}[Exponential decay]
Certain optimization problems have a graph structure:
in that case the primal-dual variable can be reordered
following a set of nodes $\mathcal{V}$ such that $w = \{w_i \}_{i \in \mathcal{V}}$.
As a consequence, the Jacobian matrix $M$ appearing in \eqref{eq:fiaccoderiv}
can also be reordered following the graph structure: in the language of linear algebra,
there exists a permutation matrix $P$ such that the matrix $P^\top M P$ is
\emph{banded}. It is known that terms in the inverse of a banded matrix decay
exponentially fast away from the diagonal~\cite{demko1984decay}.
Then, assuming SSOSC, SCS and LICQ the solution of \eqref{eq:fiaccoderiv}
satisfies the \emph{exponential decay of sensitivity} property~\cite{shin2022exponential}, that is,
for two parameters $p, p'$ close enough,
\begin{equation}
  \label{eq:exponentialdecay}
  \| w_i(p) - w(p') \| \leq \sum_{j \in \mathcal{V}} C_{ij} \| p_j - p_j' \| \; ,
  \quad \forall i \in \mathcal{V} \; ,
\end{equation}
with $w_i(p)$ the primal-dual nodal solution at node $i \in \mathcal{V}$
and $C_{ij} = \Gamma \rho^{d_G(i, j)}$ (with $\Gamma > 0$
and $\rho > 0$ two positive constants) and $d_G(i, j)$ the graph distance between node $i$ and node $j$.
The Equation~\eqref{eq:exponentialdecay} has important consequence for computing
the sensitivities in a decentralized fashion~\cite{valenzuela2024decentralized} or building a preconditionner
to solve \eqref{eq:fiaccoderiv} with an iterative method~\cite{shin2020overlapping}.
\end{remark}

\subsection{Robinson's generalized equations}
\label{sec:solutionsensitivity:robinson}
We emphasize that the strict complementarity condition is unlikely to hold on practical optimization
problems, where constraints can change locally from active to inactive.
For that reason, it has been investigated how to generalize Fiacco's result
without the strict complementarity assumption.
In \cite{robinson1979generalized}, Robinson proposed a powerful
framework to characterize the solution of nonlinear program~\eqref{eq:problem},
by interpreting the KKT conditions~\eqref{eq:kkt} as a \emph{generalized
equation} $F(w, p) + N_C(w) \ni 0$ (see Appendix~\ref{sec:robinson}
for a brief recall on generalized equations).

\begin{proposition}[KKT system as generalized equation]
  Let $w = (x, y, z) \in \mathbb{R}^{n} \times \mathbb{R}^{m_e} \times
  \mathbb{R}^{m_i}$ be a primal-dual variable.
  A primal-dual solution $w$ of the KKT system~\eqref{eq:kkt} is also solution of the generalized equation
  \begin{equation}
    \label{eq:kktgeneq}
    F(w, p) + N_C(w) \ni 0 \; ,
  \end{equation}
  with $F(w, p) := \begin{bmatrix}
    \phantom{-}\nabla_x \lagrangian(x, y, z, p) \\
      -\nabla_y \lagrangian(x, y, z, p) \\
      -\nabla_z \lagrangian(x, y, z, p)
    \end{bmatrix}$ and
    $C = \mathbb{R}^{n} \times \mathbb{R}^{m_e} \times \mathbb{R}^{m_i}_+$.
\end{proposition}
When writing the KKT system as a generalized equation, the nonsmoothness
arising from the active set changes is encoded implicitly by the normal
cone $N_C(x)$, $F$ being here a smooth functional. If the generalized equation \eqref{eq:kktgeneq}
is strongly regular, then the sensitivities follow from
Robinson's Theorem (Theorem~\ref{thm:robinsonift} in the appendix).
Strong regularity is implied by SSOSC and LICQ, allowing to drop the SCS condition
(however, SSOSC is necessary to guarantee that the perturbed solution $x(p)$ remains unique without SCS).

\begin{proposition}[Theorem 4.1, \cite{robinson1980strongly}]
  \label{thm:robinsonkkt}
  Let $w^\star = (x^\star, y^\star, z^\star)$ a primal-dual solution of the generalized
  equation~\eqref{eq:kktgeneq}. If SSOSC and LICQ hold at $w^\star$, then
  \eqref{eq:kktgeneq} is strongly regular at $w^\star$.
\end{proposition}
The strong regularity given by Proposition~\ref{thm:robinsonkkt} implies
\begin{itemize}
  \item[(a)] The solution $x^\star$ is an isolated local minimum of \nlp\
    and the associated Lagrange multipliers $(y^\star, z^\star)$ are unique;
  \item[(b)] For $p$ near $p^\star$, there exists a unique Lipschitz
    continuous function $w(p) = (x(p), y(p), z(p))$ satisfying
    KKT and SOSC for \nlp, with $w(p^\star) = w^\star$, and $x(p)$
    being a unique solution of \nlp;
  \item[(c)] LICQ holds at $x(p)$ for $p$ near $p^\star$.
\end{itemize}
Hence, strict complementarity slackness is not required to ensure
the primal-dual solution $w(p)$ is locally Lipschitz continuous.
In addition, we can prove that the optimal solution
is also directionally differentiable in a direction $h \in \mathbb{R}^{\ell}$.

\begin{proposition}[Theorem 2, \cite{jittorntrum1984solution}]
  \label{thm:jittorntrum1984}
  Suppose that at a KKT stationary solution $x^\star$
  of NLP$(p^\star)$, LICQ and SSOSC hold. Then
  \begin{itemize}
    \item there exists a local
    unique continuous primal-dual solution $w(\cdot)$ which is directionally differentiable
    in every direction $h\in \mathbb{R}^{\ell}$;
  \item for a given direction $h \in \mathbb{R}^{\ell}$, the directional derivative
    $w'(p; h)$ is the primal-dual solution of the QP problem defined as
    \begin{equation}
      \label{eq:qpjittorntrum}
      \begin{aligned}
        \min_{d}\; & \frac{1}{2} d^\top (\nabla^2_{x x}\lagrangian)  \, d + h^\top (\nabla^2_{x p}\lagrangian) \,d \\
        \text{s.t.} \quad & \nabla_x g_i(x, p)^\top d + \nabla_p g_i(x, p)^\top h = 0 & (\forall i \in [m_e]) \,,\\
                          & h_j(x, p) + \nabla_x h_j(x, p)^\top d + \nabla_p h_j(x, p)^\top h \leq 0 & (\forall j \in [m_i])\,,
      \end{aligned}
    \end{equation}
    with $\nabla_{x x}^2 \lagrangian = \nabla_{x x}^2 \lagrangian(x, y, z, p)$, $\nabla_{x p}^2 \lagrangian  = \nabla_{xp}^2 \lagrangian(x, y, z, p)$.
  \end{itemize}
\end{proposition}
The QP is parameterized implicitly by the dual multipliers $(y, z)$ (appearing
in the derivatives of the Lagrangian $\nabla_{x x}^2 \lagrangian$ and $\nabla_{x p}^2 \lagrangian$). LICQ guarantees that the multipliers
are unique at the solution, leading to a non-ambiguous definition of the
directional derivative~\eqref{eq:qpjittorntrum}.
Under LICQ and SSOSC, the directional derivative exists for every direction
$h \in \mathbb{R}^{\ell}$. However, we cannot guarantee that the directional derivative
$w'(p; \cdot)$ is continuous unless SCS hold.
Note that under SCS, the QP~\eqref{eq:qpjittorntrum}
rewrites with only equality constraints and
becomes equivalent to the linear system~\eqref{eq:fiaccoderiv}.

\subsection{Lexicographic derivative as an alternative to Robinson's generalized equation}
\label{sec:solutionsensitivity:lexicographic}
Robinson's generalized equation is a powerful framework to encode
nonsmooth equations. Several alternatives have been developed in the recent years,
one of the most promising method being based on the notion of lexicographic derivative \cite{nesterov2005lexicographic},
which has been introduced has a practical
differentiation method for nonsmooth functions, such as $PC^1$ functions (see Appendix~\ref{sec:deriv:pc1}).
The KKT conditions~\eqref{eq:kkt} rewrites equivalently as a system of nonsmooth equations:
\begin{equation}
  \label{eq:kkt:nonsmooth}
  \Phi(x, y, z, p) :=
  \begin{bmatrix}
    \nabla_x L(x, y, z, p)  \\
    g(x, p) \\
    \min\{-h(x, p), z \}
  \end{bmatrix}
  = 0  \; .
\end{equation}
It is well known that the nonsmooth KKT system \eqref{eq:kkt:nonsmooth} can be solved using a semi-smooth Newton method.
At the solution, the sensitivity analysis of the nonsmooth system can be carried out
using lexicographic derivatives~\cite{barton2018computationally}.

The function $\min$ is $PC^1$ in the sense of Scholtes \cite{scholtes2012introduction},
so we deduce that the mapping $\Phi$ is itself $PC^1$ on its domain.
The notion of coherently oriented structure extends the Implicit
Function Theorem~\ref{thm:ift} to $PC^1$ functions.

\begin{proposition}[Theorem 5.1 \cite{stechlinski2018generalized}]
  Let $(x^\star, y^\star, z^\star)$ be a solution of KKT.
  If $\Phi$ is \emph{coherently oriented} with respect to $(x, y, z)$
  at $(x^\star, y^\star, z^\star)$ then there exists a neighborhood
  $U$ of $p$ and $V$ of $x^\star$ and a $PC^1$ mapping $w: U \to V$ such that for each $p \in U$
  $w(p)$ is the unique solution of~\eqref{eq:kkt:nonsmooth}. Moreover, for any $k \in \mathbb{N}$
  and any $R \in \mathbb{R}^{p \times k}$ the LD-derivatives
  $x'(p; R)$, $y'(p, R)$ and $z'(p, R)$ are the unique
  solution $X$, $Y$ and $Z = [Z^+~Z^0~Z^-]$ of the following nonsmooth
  equation system
  \begin{equation}
    \label{eq:sensitivity:nonsmooth}
    \begin{aligned}
     \begin{bmatrix}
      \nabla^2_{x x} \lagrangian & G_x^\top & (H_x^+)^\top \\
      G_x & 0 & 0 \\
      H_x^+ & 0 & 0
    \end{bmatrix}
    \begin{bmatrix}
      X \\ Y \\ Z^+
    \end{bmatrix}
    = -
    \begin{bmatrix}
      \nabla^2_{xp} \lagrangian \\
      G_p \\
      H_p^+
    \end{bmatrix}
    R \; , \\
     \textbf{LMmin}(-H_p^0 R - H_x^0 X, Z^0) = 0 \; , \\
     Z^- = 0 \;,
    \end{aligned}
  \end{equation}
  where $\textbf{LMmin}$ is the lexicographic matrix minimum function
  defined in \eqref{eq:deriv:lmmin}, and
  \[
    G_{(x,p)} = J_{(x, p)} g(x, p) \; , \quad
    H_{(x,p)}^+ = J_{(x, p)} h_{\mathcal{I}_p^+}(x, p) \; , \quad
    H_{(x,p)}^0 = J_{(x, p)} h_{\mathcal{I}_p^0}(x, p) \; .
  \]
\end{proposition}

Without surprise, ensuring complete coherently orientation
is equivalent to assuming LICQ and SSOSC hold.
\begin{proposition}
  Let $(x^\star, y^\star, z^\star)$ be a KKT stationary point.
  Suppose that LICQ and SSOSC holds at $(x^\star, y^\star, z^\star)$.
  Then $\Phi$ is coherently oriented with respect to $(x, y, z)$.
\end{proposition}

Solving the nonsmooth system~\eqref{eq:sensitivity:nonsmooth} is generally non trivial,
and falls back to a nonsmooth Newton method or an active set algorithm (cycling through linear equation
system solves parameterized by the current working active set). However, it has been proved
in \cite{stechlinski2019generalized} that evaluating the LD derivatives can be rewritten
as the solution a hierarchy of $k$ QP problems
$QP_{(1)}, QP_{(2)}, \cdots, QP_{(k)}$, with $QP_{(j)}(r_{(j)})$
defined for $j=1, \cdots, k$ as
\begin{equation*}
  \begin{aligned}
  \min_{d\in\mathbb{R}^{n}} \; & \frac 12 d^\top \nabla_{ x x}^2 \lagrangian \, d + r_{(j)}^\top \nabla^2_{px} \lagrangian \, d \\
  \text{subject to} \quad & \nabla_x g_i(x, p)^\top d + \nabla_p g_i(x, p)^\top r_{(j)} = 0 & \forall i \in [m_e] && \leftarrow \text{multiplier}~\alpha \\
                          &  \nabla_x h_i(x, p)^\top d + \nabla_p h_i(x, p)^\top r_{(j)}= 0 & \forall i \in \mathcal{A}^+_{(j-1)} && \leftarrow \text{multiplier}~\beta \\
                          &  \nabla_x h_i(x, p)^\top d + \nabla_p h_i(x, p)^\top r_{(j)}\leq 0 & \forall i \in \mathcal{A}^0_{(j-1)} && \leftarrow \text{multiplier}~\gamma
  \end{aligned}
\end{equation*}
We note $\alpha_{(j)} \in \mathbb{R}^{m_e}$, $\beta_{(j)} \in \mathbb{R}^{m^+_i}$,
$\gamma_{(j)} \in \mathbb{R}^{m^0_i}$
the unique multipliers associated resp. to the equality constraints,
the strongly active constraints and the weakly active constraints.
The sets of strongly and weakly active constraints change between two consecutive
indexes $j$ using the information brought by the latest multiplier $\gamma_{(i)}$:
\begin{equation}
  \begin{aligned}
  & \mathcal{A}_{(j)}^0  = \{ i \in \mathcal{A}_{(j-1)}^0 \;:\; \nabla_x h_i(x, p)^\top d + \nabla_p h_i(x, p)^\top p = 0 ~\text{and}~\gamma_i = 0 \} \; , \\
  & \mathcal{A}_{(j)}^+  = \{ i \in \mathcal{A}_{(j-1)}^0 \;:\; \nabla_x h_i(x, p)^\top d + \nabla_p h_i(x, p)^\top p = 0 ~\text{and}~\gamma_i > 0 \} \cup \mathcal{A}_{(j-1)}^+ \;,
  \end{aligned}
\end{equation}
with $\mathcal{A}_{(0)}^0 = \activeset_p^0(x; z)$ and $\mathcal{A}_{(0)}^+ = \activeset_p^+(x; z)$.

\begin{proposition}[Theorem 3.1 \cite{stechlinski2019generalized}]
  Let $(x^\star, y^\star, z^\star)$ be a solution of KKT.
  Suppose that LICQ and SSOSC hold. Then, for any $k \in \mathbb{N}$
  and $R \in \mathbb{R}^{\ell \times k}$, the LD-derivatives $x'(p, R)$,
  $y'(p, R)$ and $z'(p, R)$ in the direction $R$ are given as
  \begin{equation}
    \begin{aligned}
      & x'(p, R) = \big[
        d_{(1)}(r_{(1)}) ~
        d_{(2)}(r_{(2)}) ~
        \cdots~
      d_{(k)}(r_{(k)}) \big]  \; , \\
      & y'(p, R) = \big[
        \alpha_{(1)}(r_{(1)}) ~
        \alpha_{(2)}(r_{(2)}) ~
        \cdots~
      \alpha_{(k)}(r_{(k)}) \big]  \; , \\
      & z'(p, R) = \big[
        \mu_{(1)}(r_{(1)}) ~
        \mu_{(2)}(r_{(2)}) ~
        \cdots~
      \mu_{(k)}(r_{(k)}) \big]  \; ,
    \end{aligned}
  \end{equation}
  with, for $j = 1, \cdots, k$,
  \begin{equation}
    \begin{aligned}
      & \mu_{(j)}(r_{(j)})_{\mathcal{A}^+_{(j)}} = \beta_{(j)}(r_{(j)})\;, \\
      & \mu_{(j)}(r_{(j)})_{\mathcal{A}^0_{(j)}} = \gamma_{(j)}(r_{(j)}) \;,\\
      & \mu_{(j)}(r_{(j)})_{\mathcal{A}^-_{(j)}} = 0 \; .
    \end{aligned}
  \end{equation}
\end{proposition}

\subsection{Towards degeneracy: dropping the LICQ condition}
\label{sec:solutionsensitivity:degenerate}
We have shown that the sensitivities remain well-defined if we drop the SCS
condition required in Theorem~\ref{thm:fiacco}. What happens if we suppose
in addition that the multipliers $(y, z)$ are non unique?
It turns out that in that case we retain the directional differentiability of the primal solution $x(p)$.
However, the definition of QP~\eqref{eq:qpjittorntrum}
becomes ambiguous: which multipliers should we consider when evaluating
the Hessian of the Lagrangian?

Kojima was the first to prove the continuity of the optimal solution
assuming only MFCQ and GSSOSC, and proved in addition that they are the weakest conditions
under which the perturbed solution is locally unique~\cite{kojima1980strongly}.

\begin{theorem}[Theorem 7.2, \cite{kojima1980strongly}]
  Suppose that at a KKT stationary solution $x^\star$
  of NLP$(p^\star)$, MFCQ and GSSOSC hold. Then,
  there are open neighborhoods
  $U$ of $p^\star$ and $V$ of $x^\star$ and a function $x: U \to V$ such
  that $x(\cdot)$ is continuous. The function $x(p)$ is the unique local solution
  of NLP($p$) in $V$, and MFCQ holds at $x(p)$.
\end{theorem}
The LICQ condition in Proposition~\ref{thm:jittorntrum1984} can be relaxed
to SMFCQ, which also implies the multipliers $(y^\star, z^\star)$ at the solution are unique
(Proposition~\ref{prop:kyparisis1985}).
\begin{proposition}[Theorem 4.2, \cite{shapiro1985second}]
  Suppose that at a solution $x^\star$
  of NLP$(p^\star)$, SMFCQ and SSOSC hold. Then,
  $x(\cdot)$ is directionally differentiable at $p$
  for every $h \in \mathbb{R}^{\ell}$, and the directional derivative $x'(p^\star, h)$
  is the unique solution of the QP problem
  \begin{equation}
    \label{eq:qpshapiro}
    \min_{d}\;  \frac{1}{2} d^\top (\nabla^2_{x x}\lagrangian)  \, d + h^\top (\nabla^2_{x p} \lagrangian) \, d  \quad \text{s.t.} \quad  d \in \criticalcone_{p^\star}(x^\star, z^\star; h) \;,
  \end{equation}
  with $\criticalcone_{p^\star}$ the critical cone defined in \eqref{eq:criticalcone}
  and $\nabla_{x x}^2 \lagrangian = \nabla_{x x}^2 \lagrangian(x^\star, y^\star, z^\star, p^\star)$,
  $\nabla_{x p}^2 \lagrangian = \nabla_{x p}^2 \lagrangian(x^\star, y^\star, z^\star, p^\star)$.
\end{proposition}

Note that we are no longer able to evaluate the directional derivative
of the dual solution in \eqref{eq:qpshapiro}, in contrast with the QP~\eqref{eq:qpjittorntrum}.
Kyparisis~\cite{kyparisis1990sensitivity} proved that we can relax SMFCQ by CRCQ
to get the existence of the directional derivative $x'(p, h)$ in the degenerate case
where the multipliers $(y, z)$ are non unique: it suffices
to look at the multipliers at the extreme points $\extremeset_{p^\star}(x^\star)$.

\begin{proposition}[Theorem 2.2, \cite{kyparisis1990sensitivity}]
  \label{thm:sensitivity:kyparisis}
  Suppose that at a KKT stationary solution $x^\star$
  of NLP$(p^\star)$, MFCQ, CRCQ and GSSOSC hold. Then,
  $x(\cdot)$ is directionally differentiable at $p^\star$
  for every $h \in \mathbb{R}^{\ell}$, and the directional derivative $x'(p^\star, h)$
  is the unique solution of the QP problem, defined for a given $(y, z) \in \extremeset_{p^\star}(x^\star)$,
  \begin{equation}
    \label{eq:qpkyparisis}
    \min_{d}\;  \frac{1}{2} d^\top (\nabla_{x x}^2 \lagrangian) \, d + h^\top (\nabla_{x p}^2 \lagrangian) d  \quad \text{s.t.} \quad  d \in \criticalcone_{p^\star}(x^\star, z; h) \;,
  \end{equation}
  with $\criticalcone_{p^\star}$ the critical cone defined in \eqref{eq:criticalcone}
  and $\nabla_{x x}^2 \lagrangian = \nabla_{x x}^2 \lagrangian(x^\star, y, z, p^\star)$,
  $\nabla_{x p}^2 \lagrangian = \nabla_{x p}^2 \lagrangian(x^\star, y, z, p^\star)$.
\end{proposition}

Finally, Proposition~\ref{thm:sensitivity:kyparisis} has been extended
in \cite{ralph1995directional} to derive a practical way to evaluate the sensitivities for a
degenerate nonlinear program. The idea is to look only at multipliers in $\extremeset_{p^\star}(x^\star)$
that are also solution of the linear program
\begin{equation}
  \label{eq:dempe:LP}
  S_{p^\star}(x^\star; h) := \left\{ \begin{aligned}
  \argmax_{y, z} \; &
    y^\top \big(J_p g(x^\star, p^\star)\big)\, h +
    z^\top \big(J_p h(x^\star, p^\star)\big)\, h \\
  \text{s.t.} \quad & (y, z) \in \multiplierset_{p^\star}(x^\star) \;.
  \end{aligned}
  \right.
\end{equation}
The solution $(y, z)$ of \eqref{eq:dempe:LP} determines the sets of strongly active
$\activeset_p^+(x^\star; z)$ and weakly active constraints $\activeset_p^0(x^\star; z)$ in the QP~\eqref{eq:qpkyparisis}.
As such, the set $S_{p^\star}(x^\star; h)$ is intimately related to the critical set $\criticalcone_p(x^\star, z; h)$.
Indeed, observe that by definition $(y, z) \in \multiplierset_p(x)$ if
\begin{equation}
  \label{eq:dempe:cond2}
     \nabla_x f(x^\star, p^\star) + \sum_{i=1}^{m_e} y_i \nabla_x g_i(x^\star, p^\star) + \sum_{i \in \activeset_p(x^\star)}z_i \nabla_x h_i(x^\star, p^\star) = 0
     \;, \quad  z_{\mathcal{I}} \geq 0  \;,
\end{equation}
where we have removed the inactive inequality from \eqref{eq:dempe:cond2}.
By noting $d \in \mathbb{R}^n$ the multiplier associated to the stationary condition,
and $\lambda \in \mathbb{R}^{|\activeset_p(x^\star)|}$ the multiplier associated to
the constraints $z_\mathcal{I} \geq 0$, the KKT conditions of \eqref{eq:dempe:LP} write
\begin{equation*}
  \begin{aligned}
     & \nabla_x g_i(x^\star, p^\star)^\top d + \nabla_p g_i(x^\star, p^\star)^\top h = 0 & \forall i \in [m_e]\;, \\
     & \nabla_x h_i(x^\star, p^\star)^\top d + \nabla_p h_i(x^\star, p^\star)^\top h + \lambda = 0  & \forall i \in \activeset_p(x^\star)\;, \\
     & \nabla_x f(x^\star, p^\star) + \sum_{i=1}^{m_e} y_i \nabla_x g_i(x^\star, p^\star) + \sum_{i \in \activeset_p(x^\star)}z_i \nabla_x h_i(x^\star, p^\star) = 0 \;,\\
     & 0 \leq \lambda \perp z_\mathcal{I} \geq 0 \; .
  \end{aligned}
\end{equation*}
By removing the multiplier $\lambda$ and taking into account the disjunction
in the complementarity constraint $0 \leq \lambda \perp z_\mathcal{I} \geq 0$, the two first equations rewrite equivalently as
\begin{equation}
  \label{eq:dempe:cond3}
  \begin{aligned}
     & \nabla_x g_i(x^\star, p^\star)^\top d + \nabla_p g_i(x^\star, p^\star)^\top h = 0 & \forall i \in [m_e]\;, \\
    & \nabla_x h_i(x^\star, p^\star)^\top d + \nabla_p h_i(x^\star, p^\star)^\top h \leq 0 & \forall i \in \activeset_p^0(x^\star; z) \;, \\
    & \nabla_x h_i(x^\star, p^\star)^\top d + \nabla_p h_i(x^\star, p^\star)^\top h = 0 & \forall i \in \activeset_p^+(x^\star; z) \; , \\
  \end{aligned}
\end{equation}
which describes exactly the critical set $\criticalcone_p(x^\star, z; h)$.
Hence, by theory of linear programming, the problem~\eqref{eq:dempe:LP}
has an optimal solution if and only if the previous KKT system is consistent,
meaning that $\criticalcone_p(x^\star,z;h)$ is nonempty.
This is formalized in the following lemma.
\begin{lemma}[Lemma 2.2, \cite{dempe1993directional}]
  The critical cone $\criticalcone_p(x^\star, z; h)$ is nonempty if and only if
  $(y, z) \in S_p(x^\star; h)$.
\end{lemma}

We can go one step further in the interpretation of the LP \eqref{eq:dempe:LP}.
Note that its dual writes
\begin{equation}
  \begin{aligned}
    \min_{d \in \mathbb{R}^n} \; & \nabla_x f(x^\star, p^\star)^\top d \; , \\
    \text{s.t.} \quad & \nabla_x g_i(x^\star, p^\star)^\top d + \nabla_p g_i(x^\star, p^\star)^\top h = 0 \; , \; \forall i \in [m_e] \; ,\\
      & \nabla_x h_i(x^\star, p^\star)^\top d + \nabla_p h_i(x^\star, p^\star)^\top h \leq 0 \; , \; \forall i \in \activeset_p(x^\star) \; ,\\
  \end{aligned}
\end{equation}
with the multiplier $y \in \mathbb{R}^{m_e}$ being associated to the equality
constraints and the multiplier $z \in \mathbb{R}^{m_i}$ being associated
to the inequality constraints. As a consequence, the multipliers $(y, z)$
are the marginal costs encoding the objective's change after a perturbation $h$,
and satisfy
\begin{equation}
  \nabla_x f(x^\star, p^\star)^\top d + \sum_{i=1}^{m_e} y_i \nabla_p g_i(x^\star, p^\star)^\top h
  + \sum_{i \in \activeset_p(x^\star)}^{m_e} z_i \nabla_p h_i(x^\star, p^\star)^\top h = 0 \; .
\end{equation}
The nature of the LP \eqref{eq:dempe:LP} has been further studied
in \cite[Lemma 4.1]{izmailov2005note}. Izmailov and Solodov noted that
the solution of \eqref{eq:dempe:LP} is likely to be unique if
the objective of \eqref{eq:dempe:LP} is non-constant. Furthermore, if SMFCQ
is not satisfied the set of multipliers $\multiplierset_p(x^\star)$ does not
reduce to a singleton. Then, the LP \eqref{eq:dempe:LP} returns
a solution $(y, z) \in \extremeset_p(x^\star)$ that does not satisfy strict-complementarity: $\activeset_p^0(x^\star, z) \neq \emptyset$.

Ralph and Dempe showed in \cite{ralph1995directional} that the directional derivative is unique
even when the solution of the LP \eqref{eq:dempe:LP} is not.
They also prove that the local solution $x(\cdot)$ is piecewise differential
in the sense of Scholtes~\cite{scholtes2012introduction}.
\begin{theorem}[Theorem 2, \cite{ralph1995directional}]
  \label{thm:sensitivity:ralphdempe}
  Suppose that at a KKT stationary solution $x^\star$
  of NLP$(p^\star)$, MFCQ, CRCQ and GSSOSC hold. Then,
  \begin{itemize}
    \item There is open neighborhoods $U$ of $p^\star$ and $V$ of $x^\star$
      and a function $x:U \to V$  such that $x(p)$ is the unique local
      solution of \nlp\ in $V$ and MFCQ holds at $x(p)$;
    \item $x(\cdot)$ is a PC$^1$ function (hence locally Lipschitz
      and B-differentiable);
    \item For each $p \in U$, the directional derivative $x'(p; \cdot)$ is a piecewise linear
      function, and for $h \in \mathbb{R}^{\ell}$
      and $(y, z) \in S_p(x^\star; h)$, $x'(p; h)$ is the unique solution of \eqref{eq:qpkyparisis}.
  \end{itemize}
\end{theorem}
It turns out that the notion of $PC^{1}$ functions is adapted to characterize
the solution of a parametric nonlinear program.

\section{Sensitivity of the optimal value function}
\label{sec:objectivesensitivity}
We now focus on the derivative of the value function $\valuefunction(p)$.
In \S\ref{sec:objectivesensitivity:nonsmooth}, we recall classical results
focusing on the directional differentiability of the value function $\valuefunction$.
The differentiable case is studied in \S\ref{sec:objectivesensitivity:smooth}.

\subsection{Danskin lemma and its extensions}
\label{sec:objectivesensitivity:nonsmooth}
The sensitivity analysis of the optimal value function
traces back to the seminal work of Danskin, who studied the sensitivity
for a simplified problem where the admissible set $X(p)$ is a
fixed subset $X \subset \mathbb{R}^{n}$, that is,
\begin{equation}
  \label{eq:danskinpb}
  \valuefunction(p) = \min_{x \in X} \; f(x, p) \; .
\end{equation}
\begin{theorem}[Danskin's Lemma~\cite{danskin2012theory}]
  Suppose $X$ is nonempty and compact, and $f(x, \cdot)$ is continuously differentiable
  at $p \in \mathbb{R}^{\ell}$.
  Then $\valuefunction$ is locally Lipschitz near $p$,
  and directionally differentiable in every
  direction $h \in \mathbb{R}^{\ell}$, with
  \begin{equation}
    \valuefunction'(p; h) = \min_{x \in \solutionset(p)} \; \nabla_p f(x, p)^\top h \; ,
  \end{equation}
  with $\solutionset(p)$ the solution set defined in \eqref{eq:solutionset}.
\end{theorem}
Danskin's theorem has been extended by Hogan to the case
where $X(p)$ is defined as $X(p) = \{ x \in \mathbb{R}^{n} \;:\; x \in M \;,\; h(x, p) \leq 0 \}$,
for a fixed set $M$.
\begin{theorem}[\cite{hogan1973point}]
  Suppose $M$ is a closed convex set, $f(\cdot, p)$ and $h_i(\cdot, p)$ are convex on $M$ for all $p$, and
  continuously differentiable on $M \times N$, with $N$ a neighborhood of $p$.
  If (i) the solution set $\solutionset(p)$ is nonempty and bounded, (ii) $\valuefunction(p)$ is finite,
  (iii) there is a point $\hat{x} \in M$ such that $h(\hat{x}, p) < 0$ (Slater),
  then $\valuefunction$ is directionally differentiable, and for all $h\in \mathbb{R}^{\ell}$,
  \begin{equation}
  \valuefunction'(p; h) = \min_{x \in \solutionset(p)} \max_{z \in \multiplierset_p(x)} \nabla_p \lagrangian(x, z, p)^\top h
    \;,
  \end{equation}
  with $\multiplierset_p(x)$ the multiplier solution set.
\end{theorem}

The extension to the generic case
associated to $X(p) = \{x \in \mathbb{R}^{n} \; | \; g(x, p) = 0 \, , \, h(x, p) \leq 0  \}$ was proved by Gauvin and Dubeau~\cite{gauvin1982differential}
and Fiacco~\cite{fiacco1982optimal}.
The following theorem bounds the lower and upper Dini derivatives (see
Appendix \ref{sec:deriv:nonsmooth})
of the value function: in general the bounds are sharp.
\begin{theorem}[Theorem 2.3.4 \cite{fiacco1983introduction}]
  Suppose that $X(p)$ is nonempty and uniformly compact near $p$
  and MFCQ holds at each $x \in \solutionset(p)$. Then $\valuefunction$ is locally Lipschitz
  near $p$, and for any $h \in \mathbb{R}^{\ell}$,
  \begin{equation}
    \begin{aligned}
    \inf_{x \in \solutionset(p)} \min_{(y, z) \in \multiplierset_p(x)} \nabla_p \lagrangian(x, y, z, p)^\top h
    & \leq \lim_{t \downarrow 0} \inf \frac{1}{t}\big(\valuefunction(p + th) - \valuefunction(p)\big) \\
    & \leq \lim_{t \downarrow 0} \sup \frac{1}{t}\big(\valuefunction(p + th) - \valuefunction(p)\big) \\
    & \leq \inf_{x \in \solutionset(p)} \max_{(y, z) \in \multiplierset_p(x)} \nabla_p \lagrangian(x, y, z, p)^\top h \;.
    \end{aligned}
  \end{equation}
\end{theorem}

In the convex case the multiplier set $\multiplierset_p(x)$ is
independent on the solution $x \in \solutionset(p)$.
\begin{proposition}[Corollary 2.3.8 \cite{fiacco1983introduction}]
  Let $f$ and $h$ be convex functions in $x$, $g$ be affine in $x$, such that
  $f,g, h$ are all jointly $C^1$ in $(x, p)$. Then, if $X(p)$ is nonempty
  and uniformly compact near $p$ and MFCQ holds at each $x \in \solutionset(p)$,
  then $\valuefunction$ is directionally differentiable at $p$ and
  \begin{equation}
    \valuefunction'(p; h) = \min_{x \in \solutionset(p)} \max_{(y, z) \in \multiplierset_p(x)}
    \nabla_p \lagrangian(x, y, z, p)^\top h \; .
  \end{equation}
\end{proposition}

Note that if in addition SMCFQ holds, then $\multiplierset_p(x)$ reduces to a singleton
$(y^\star, z^\star)$ and the directional derivative simplifies to
\begin{equation}
  \valuefunction'(p; h) = \min_{x \in \solutionset(p)}
  \nabla_p \lagrangian(x, y^\star, z^\star, p)^\top h \; .
\end{equation}
If instead we suppose GSSOSC holds, the local solution $x(p)$ is a singleton
and the computation of the directional derivative also simplifies.
\begin{theorem}[Theorem 7.5 \cite{fiacco2006sensitivity}]
  Suppose GSSOSC and MFCQ hold at a stationary solution $x^\star$.
  Then near $p$ the function $\valuefunction(\cdot)$ has a finite one-sided
  directional derivative, and for every $h \in \mathbb{R}^{\ell}$,
  \begin{equation}
    \valuefunction'(p; h) = \max_{(y, z) \in \multiplierset_p(x^\star)} \nabla_p \lagrangian(x^\star, y, z, p)^\top h
    \; .
  \end{equation}
\end{theorem}

\subsection{The differentiable case}
\label{sec:objectivesensitivity:smooth}
Now, we impose more stringent constraint qualifications to recover the setting of
the Implicit Function Theorem~\ref{thm:ift}: as a consequence the value function $\valuefunction$ becomes $C^2$ differentiable.
For a given $p \in \mathbb{R}^{\ell}$,
the following results hold
not only for global solution $x \in \solutionset(p)$, but also \emph{locally} for every
KKT stationary point $w(p) = (x(p), y(p), z(p))$ of \nlp.
As such, we define the \emph{local value function} as
\begin{equation}
  \valuefunction_\ell(p) = f(x(p), p) \; .
\end{equation}

\begin{theorem}[Theorem 3.4.1 \cite{fiacco1983introduction}]
  \label{thm:valuefunction:fiacco}
  Suppose SOSC, SCS and LICQ hold at a stationary KKT solution $(x^\star, y^\star, z^\star)$.
  Then the local value function $\valuefunction_\ell$ is $C^2$ near $p$, and
  \begin{itemize}
    \item $\valuefunction_\ell(p) = \lagrangian(x^\star, y^\star, z^\star, p)$.
    \item $\nabla_p \valuefunction_\ell(p) = \nabla_p \lagrangian(x^\star, y^\star, z^\star, p)$.
    \item and also
      \begin{equation}
        \begin{aligned}
          \nabla^2_{pp} \valuefunction_\ell(p) = \nabla_p \big(\nabla_p \lagrangian(w^\star, p)\big)
                                                &= \nabla_{pp}^2 \lagrangian + \nabla_{wp}^2 \lagrangian \; J_p w(p) \;, \\
                                                &= \nabla_{pp}^2 \lagrangian - \nabla_{wp}^2 \lagrangian \, (\nabla_{ww}^2 \lagrangian)^{-1} \nabla_{pw}^2 \lagrangian \;.
        \end{aligned}
      \end{equation}
  \end{itemize}
\end{theorem}
The previous theorem gives expressions both for the gradient and the Hessian
of the optimal value function. It applies to generic nonlinear programs of the form
\nlp. It yields interesting extensions when applied to problem with certain structures.

\paragraph{Case 1: the constraints are independent of $p$}
First, we can prove that if the constraints of the nonlinear problem
are independent of the parameter $p$,
\begin{equation}
  \label{eq:problem_struct1}
  \text{NLP}_{obj}(p) := ~
  \min_{x \in \mathbb{R}^{n}} \; f(x, p) \quad \text{subject to} \quad
  \left\{
  \begin{aligned}
    g(x) = 0 \;,\\
    h(x) \leq 0 \;,
  \end{aligned}
  \right.
\end{equation}
then the gradient and the Hessian depend only on the derivatives of the objective function~$f$.

\begin{corollary}
  Suppose the conditions of Theorem~\ref{thm:valuefunction:fiacco} hold
  at a stationary solution $(x^\star, y^\star, z^\star)$ of problem \eqref{eq:problem_struct1}.
  Then, in a neighborhood
  of $p$ the local value function $\valuefunction_\ell$ is $C^2$, with
  \begin{itemize}
    \item $\valuefunction_\ell(p) = f(x^\star, p)$.
    \item $\nabla_p \valuefunction_\ell(p) = \nabla_p f(x^\star)$.
    \item $\nabla^2_{pp} \valuefunction_\ell(p) = \nabla_{pp}^2 f + (\nabla_{xp}^2 f) \, J_p x(p)$.
  \end{itemize}
\end{corollary}

\paragraph{Case 2: RHS perturbations and shadow prices}
The other corollary draws a link between Theorem~\ref{thm:valuefunction:fiacco}
and the shadow price of sensitivity: if the parameters appear only in the right-hand-side
of the constraint, that is, for $p \in \mathbb{R}^{m_e}$ and $q \in \mathbb{R}^{m_i}$,
\begin{equation}
  \label{eq:problem_struct2}
  \text{NLP}_{rhs}(p, q) := ~
  \min_{x \in \mathbb{R}^{n}} \; f(x) \quad \text{subject to} \quad
  \left\{
  \begin{aligned}
    g(x) = p \;,\\
    h(x) \leq q \;,
  \end{aligned}
  \right.
\end{equation}
then the gradient of the value function is equal to the Lagrange multipliers
found at the solution.
This is well-known result in convex analysis, the problem~\eqref{eq:problem_struct2}
being subject to the so-called \emph{canonical perturbations}.

\begin{corollary}[Optimal value function derivatives for RHS perturbations] \\
  Suppose the conditions of Theorem~\ref{thm:valuefunction:fiacco} hold
  at a stationary solution $(x^\star, y^\star, z^\star)$ of problem \eqref{eq:problem_struct2}.
  Then, in a neighborhood of $(p, q)$ the local value function $\valuefunction_\ell$ is $C^2$, with
  \begin{equation}
    \nabla \valuefunction_\ell(p, q) = \begin{bmatrix} y^\star \\ z^\star \end{bmatrix}
    \quad \text{and} \quad
    \nabla^2 \valuefunction_\ell(p^\star) =
    \begin{bmatrix}
      J_p y(p, q) & J_q y(p, q)  \\
      J_p z(p, q) & J_q z(p, q)
    \end{bmatrix} \; .
  \end{equation}
\end{corollary}

\section{Numerical algorithms for sensitivity analysis}
\label{sec:numerics}
Now the theory has been laid out, we give an overview of the existing numerical
methods for the sensitivity analysis of a nonlinear optimization program~\eqref{eq:problem}.
First, we put a special emphasis on the differentiation of conic programs in \S\ref{sec:num:conic},
before focusing on the differentiation of nonlinear programs in \S\ref{sec:num:nlp}.

\subsection{Differentiation of conic programs}
\label{sec:num:conic}
Conic programming has become an important subfield of optimization.
In machine learning, a significant number of problems are
conic, motivating a renewed interest for the sensitivity analysis
of conic-structured problems.
A conic program is a specialization of \eqref{eq:problem}, with the following
structure
\begin{equation}
  \label{eq:problem:conic}
    \min_{x \in \mathbb{R}^n} \;  c^\top x \quad \text{s.t.} \quad  A x = b \;, \quad x \in \mathcal{K} \;,
\end{equation}
where $\mathcal{K}$ is a pointed convex cone.
If the cone is the non-negative orthant $\mathcal{K} = \mathbb{R}^n_+$, the problem \eqref{eq:problem:conic} becomes a LP.
We suppose here that the problem's data are the parameters: $p = (c, A, b)$.

The KKT conditions of \eqref{eq:problem:conic} write out
\begin{equation}
  \label{eq:conic:kkt}
  A^\top y + s = c \;, \quad
  A x = b \; , \quad
  x \in K \; , \quad
  s \in K^\star \; , \quad
  x^\top s = 0  \; .
\end{equation}
As the problem~\eqref{eq:problem:conic} is convex, the KKT conditions \eqref{eq:conic:kkt}
are necessary and sufficient.

\paragraph{Homogeneous self-dual embedding}
On the contrary to what we exposed in Theorem~\ref{thm:fiacco},
the sensitivity analysis of \eqref{eq:problem:conic}
does not differentiate through the KKT conditions~\eqref{eq:conic:kkt}:
instead, it uses an equivalent form of~\eqref{eq:conic:kkt}
called the \emph{homogeneous self-dual embedding} (HSD) \cite{bolte2021nonsmooth,busseti2019solution}.
The HSD writes out as the problem
\begin{equation}
  \label{eq:HSD}
  \text{Find} \; (x, y, \tau, s, \kappa) \quad \text{s.t.} \quad
  \left\{
    \begin{aligned}
  & \begin{bmatrix}
    0 & -A^\top & c \\
    A & 0 & -b \\
    -c^\top & b^\top & 0
  \end{bmatrix}
  \begin{bmatrix}
    x \\ y \\ \tau
  \end{bmatrix}
  =
  \begin{bmatrix}
    s \\ 0 \\ \kappa
  \end{bmatrix} \; , \\
  & x \in K\; , \; s \in K^\star \\
  & (\tau, \kappa) \geq 0 \quad  \text{and} \; x^\top s + \tau\kappa = 0 \; .
    \end{aligned}
    \right.
\end{equation}
If $(x, y, \tau, s, \kappa)$ is solution of \eqref{eq:HSD}
with $\tau > 0$ and $\kappa = 0$, then $(x/\tau, y/\tau, s/\tau)$
satisfies the KKT conditions~\eqref{eq:conic:kkt} and is solution of
the conic program~\eqref{eq:problem:conic}.
In general, the parameters $\tau$ and $\kappa$ characterize the primal and dual
infeasibility of the conic problem~\eqref{eq:problem:conic}.

We define the skew-symmetric matrix
\begin{equation}
  Q = \begin{bmatrix}
    0 & -A^\top & c\\
    A^\top & 0 & -b \\
    -c^\top & b^\top & 0
  \end{bmatrix}
  \;,
\end{equation}
Defining the new cone $\mathcal{C} = K \times \mathbb{R}^m \times \mathbb{R}_+$ and its dual
$\mathcal{C}^\ast = K^\star \times \{0 \}^n  \times \mathbb{R}_+$,
the problem \eqref{eq:HSD} becomes equivalent to
\begin{equation}
  \label{eq:conic:hsde}
  \text{Find} \; (u, v) \in \mathcal{C} \times \mathcal{C}^\star
  \quad \text{s.t} \quad Qu = v \; , \quad u^\top v = 0 \; .
\end{equation}

\paragraph{Residual map}
Let $z \in \mathbb{R}^{n+m+1}$. Using Moreau's decomposition theorem,
we know that $z = u - v$ with $(u, v) \in \mathcal{C} \times \mathcal{C}^\star$ and $u^\top v = 0$
if and only if $u = P_{\mathcal{C}}(z)$ and $v = -P_{\mathcal{C}^o}(z)$.

As a consequence, the HSD~\eqref{eq:conic:hsde} rewrites more compactly as
finding $z \in \mathbb{R}^{n+m+1}$ such that
\begin{equation}
  -P_{\mathcal{C}^o}(z) = Q P_{\mathcal{C}}(z)  \; .
\end{equation}
The residual map $R(z)$ returns the HSD residual:
\begin{equation}
  \label{eq:conic:residualmapping}
  R(z) = Q P_{\mathcal{C}}(z) + P_{\mathcal{C}^o}(z) = \big((Q-I)P_{\mathcal{C}} + I\big) z \; .
\end{equation}
If we decompose $z^\star$ component
by component as $z^\star = (u^\star, v^\star,w^\star)$, with $(u^\star,v^\star,w^\star) \in \mathbb{R}^n \times \mathbb{R}^m \times \mathbb{R}_+^*$,
the KKT solutions $(x^\star, y^\star, s^\star)$ of \eqref{eq:problem:conic}
is recovered as
\begin{equation}
  (x^\star, y^\star, s^\star) = \phi(z^\star) \quad \text{with} \quad
  \phi(z^\star) := \frac{1}{w^\star} \big(P_{\mathcal{C}}(z^\star), -P_{\mathcal{C}^o}(z^\star), w^\star \big)
  \; .
\end{equation}

\paragraph{Sensitivity analysis}
We now have all the elements to evaluate the sensitivity at a given
KKT primal-dual solution $w^\star = (x^\star, s^\star, y^\star)$ of \eqref{eq:conic:kkt}.
By defining $u^\star = (x^\star, y^\star, 1)$,
$v^\star = (s^\star, 0, 0)$ and $z^\star = u^\star - v^\star$, the
vector $z^\star$ is solution of the residual map: $R(z^\star) = 0$.
If the projection operator $P_{\mathcal{C}}$ is differentiable,
then $R$ is itself differentiable with $J_z R(z) = (Q+I)J_z P_{\mathcal{K}}(z) - I$.
In that case, by applying the implicit function theorem \ref{thm:ift} to $R(z) = 0$, we get the Jacobian
\begin{equation}
  \label{eq:cvx:ift}
  J_p z^\star = - (J_z R(z^\star))^{-1} J_p R(z^\star) \;.
\end{equation}
We recover the sensitivity of the primal-dual solution by noting
that $w^\star =  \phi(z^\star)$: after applying
the chain-rule, we get $J_p w^\star = J_z \phi \cdot J_p z^\star$.

Note that in general the projection operator $P_{K}$
is not differentiable. However, as recalled
in \cite{ali2017semismooth}, if $K$ encodes the nonnegative orthant,
the second-order cone or the positive definite cone then
$P_{K}$ is strongly semismooth (see Appendix~\ref{sec:deriv:semismooth}), hence differentiable almost everywhere.
As a result, the residual mapping~\eqref{eq:conic:residualmapping} is a strongly semismooth function:
the equation~\eqref{eq:cvx:ift} is well-defined if we use a version of the implicit
function theorem adapted for the strongly semismooth case (and involving Clarke generalized Jacobians).
Another rigorous treatment has been made in \cite{bolte2021nonsmooth},
based on the newly introduced concept of \emph{conservative Jacobians}.

\paragraph{Implementation}
With the emergence of the differentiable programming paradigm, numerous packages
have been developed to differentiate through conic program~\eqref{eq:problem:conic}.
The approach has been implemented in {\tt cvxpylayer}~\cite{agrawal2019differentiable} and {\tt DiffOpt.jl}~\cite{sharma2022flexible}.
We emphasize that factorizing the matrix $\nabla_z R$ in \eqref{eq:cvx:ift} can be challenging,
as on the contrary to the KKT matrix the Jacobian $\nabla_z R$ is non symmetric.
For that reason, the system~\eqref{eq:cvx:ift} is solved either using a LU factorization,
or using a Krylov method (the LSQR algorithm accommodates nicely the degenerate case
where the Jacobian of the residual $\nabla_z R$ is non-invertible).

\subsection{Differentiation of nonlinear programs}
\label{sec:num:nlp}

The solution of \nlp\ has been studied by Fiacco and McCormick in their seminal
work~\cite{fiacco1990nonlinear}: the main idea is to reformulate the problem
with penalties and solve the resulting unconstrained problem using Newton method.
In that case, the sensitivity analysis uses the approximate solution returned by the algorithm
to compute the directional derivative $x'(p, h)$.

\subsubsection{Fiacco and SENSUMT}
We recall the penalty approach described originally in \cite{fiacco1980user}.
Although this approach is now slightly outdated, it introduces some fundamental
ideas still being used today.
We follow the approach described in \cite[Chapter 6]{fiacco1983introduction}.
The idea is to reformulate \nlp\ with inner and outer penalties. For a
strictly feasible point $x$ (in the sense $h(x, p) < 0$), Fiacco and McCormick define the functional
\begin{equation}
  \label{eq:fiaccopenalty}
  W(x, r, p) = f(x, p) - r \sum_{i=1}^{m_i} \log(-h_i(x, p)) + \frac{1}{2 r} \sum_{j=1}^{m_e} g_j^2(x, p) \; ,
\end{equation}
where $r > 0$ is a penalty parameter appearing in the barrier associated
to the inequality constraints and the quadratic penalty associated to the equality equations.
For a fixed penalty $r$, a stationary point $\nabla_x W\big(x(p, r), r, p\big) = 0$ is found using the
traditional Newton method.
SENSUMT~\cite{fiacco1980user} was developed as an extension of SUMT to compute
the sensitivity using the penalty function $W(\cdot)$.

Note that the gradient of the penalty $W(\cdot)$ and the gradient
of the Lagrangian \eqref{eq:lagrangian} are defined respectively as
\begin{equation}
  \begin{aligned}
    \nabla_x W(x, r, p) &= \nabla_x f(x, p) + \sum_{i=1}^{m_i} (r/h_i) \nabla_x h_i(x,p) + \sum_{j=1}^{m_e} (g_j/r) \nabla_x g_j(x, p) \; , \\
    \nabla_x \mathcal{L}(x, y, z, p) &= \nabla_x f(x, p) + \sum_{i=1}^{m_i} z_i \nabla_x h_i(x,p) + \sum_{j=1}^{m_e} y_j \nabla_x g_j(x, p) \;,
  \end{aligned}
\end{equation}
where we defined $h_i := h_i(x, p)$ and $g_j := g_j(x, p)$.
We observe that for any solution $x^\star$ of \eqref{eq:fiaccopenalty},
the multipliers $z^\star = r / h(x^\star, p)$ and $y^\star = g(x^\star, p) / r$
satisfy
the stationary condition $\nabla_x \mathcal{L}(x^\star, y^\star, z^\star, p) = 0$:
The variables $(y^\star, z^\star)$ are approximations of the optimal Lagrange multipliers,
both depending on the penalty $r$.
The penalty $r$ is here interpreted as a new parameter, and is valid candidate
to apply Theorem~\ref{thm:fiacco}.

\begin{proposition}[Theorem 6.2.1 \cite{fiacco1983introduction}]
  \label{thm:fiacco:penalty}
  Let $p^\star \in \mathbb{R}^{\ell}$.
  Suppose that a stationary solution $x^\star$ of \nlp\ LICQ, SOSC and SCS hold.
  Then in a neighborhood of $(0, p^\star)$ there exists a unique differentiable
  function $w(r, p) = (x(r, p), y(r, p), z(r, p))$ satisfying
  \begin{equation}
    \label{eq:fiaccosystem}
    \begin{aligned}
    & \nabla_x \lagrangian(x, y, z, p) = 0 \;, \\
    & g_j(x, p) = y_j r \; & \forall j \in [m_e] , \\
    & z_i h_i(x, p) = r \;, &\forall i \in [m_i] \;,
    \end{aligned}
  \end{equation}
  and such that $w(0, p^\star) = (x^\star, y^\star, z^\star)$.
  In addition, for any $(r, p)$ near $(0, p^\star)$, $x(r, p)$ is locally unique
  unconstrained local minimizing point of $W(x, r, p)$ with $h_i(x(r, p), p) < 0$
  and $\nabla^2_{x x} W(x(r, p), r, p)$ positive definite.
\end{proposition}

At the solution of \eqref{eq:fiaccopenalty}, the sensitivity of the primal $x(r, p)$ is given
by the Implicit Function Theorem \ref{thm:ift} as
\begin{equation}
  J_p x(r, p) = - (\nabla_{x x}^2 W)^{-1} \nabla_{pw} W \;,
\end{equation}
involving only the inverse of the $n \times n$ Hessian matrix $\nabla_{x x}^2 W$.
Using Proposition~\ref{thm:fiacco:penalty}, we get that $\lim_{r \to 0} J_p x(r, p) = J_p x^\star(p)$.
The dual sensitivities are recovered with~\eqref{eq:fiaccosystem} as
$z_i(r, p) = r / h_i(x(r, p), p)$ and $y_j = g_j(x(r, p), p) / r$, giving
the respective Jacobians
\begin{equation}
  \begin{aligned}
    & J_p y(r, p)= \frac{1}{r} \big(J_x g(x(r, p), p) \cdot J_p x(r, p) + J_p g(x(r, p), p)\big) \; , \\
    & J_p z(r, p) = - \frac{r}{h(x(r, p), p)} \big( J_x h(x(r, p), p) \cdot J_p x(r, p) + J_p h(x(r, p), p)\big) \;.
  \end{aligned}
\end{equation}

\subsubsection{Sensitivity in interior-point}
The SUMT algorithm described in the previous paragraph has fallen out of fashion,
as the resulting problem~\eqref{eq:fiaccopenalty} becomes highly ill-conditioned when we are approaching an optimal solution.
However, SUMT has nurtured the development of interior-point methods in the 1980s:
the primal-dual interior-point method (IPM) is now one of the most widely
used algorithm to solve nonlinear program \nlp.

The solution sensitivity of primal-dual interior point
follows from the same principle as in Theorem \ref{thm:fiacco:penalty}.
As an example, we describe the approach uses in sIpopt~\cite{pirnay2012optimal}, a package
for the sensitivity analysis of NLP developed as an extension of the
primal-dual interior-point solver Ipopt. Ipopt rewrites the
inequalities in \nlp\ as equality constraints by introducing additional slack variables,
giving a reformulated NLP problem with only bound constraints $x \geq 0$:
\begin{equation}
  \label{eq:problem:sipopt}
  \begin{aligned}
    \min_{x} \;&  f(x, p) \\
    \text{s.t.} \quad & g(x, p) = 0 \;, \quad x \geq 0 \; .
  \end{aligned}
\end{equation}
For a given barrier parameter $\mu > 0$,
the bound constraints $x \geq 0$ are penalized using a barrier term in the objective function:
\begin{equation}
  \label{eq:problem:sipoptbarrier}
  \begin{aligned}
    \min_{x} \;&  f(x, p) - \mu \sum_{i=1}^{n} \log(x_i) \\
    \text{s.t.} \quad & g(x, p) = 0  \; .
  \end{aligned}
\end{equation}
Ipopt solves a sequence of subproblems~\eqref{eq:problem:sipoptbarrier} for given barrier parameters
$\{\mu_k \}_k$, with $\mu_k \to 0$. For a fixed $\mu$, solving \eqref{eq:problem:sipoptbarrier}
with primal-dual IPM amounts to solve the system of nonlinear equations $F_\mu(x, y, z, p) = 0$, with
\begin{equation}
  F_\mu(x, y, z, p) = \begin{bmatrix}
    \nabla_x f(x, p) + J_x g(x, p)^\top y - z \\
    g(x, p) \\
    X z - \mu e
  \end{bmatrix} \; ,
\end{equation}
where we have defined $X = \text{diag}(x)$.

The following proposition is a corollary of Proposition~\ref{thm:fiacco:penalty}.
\begin{proposition}[Convergence of interior-point]
  Let $w^\star = (x^\star, y^\star, z^\star)$ be a local primal-dual stationary point of \eqref{eq:problem:sipopt}.
  Suppose that LICQ, SCS and SOSC hold at $w^\star$.
  We denote by $x(\mu_k, p)$ the solution of the barrier problem~\eqref{eq:problem:sipoptbarrier}
  at iteration $k$. Then
  \begin{itemize}
    \item There is at least one subsequence of $x(\mu_k, p)$ converging to $x^\star$
    and for every convergent subsequence, the corresponding barrier multiplier approximation
      is bounded and converges to $(y^\star, z^\star)$.
    \item There exists a unique continuously differentiable function $w(\cdot, p)$
      existing in a neighborhood of $0$ for $\mu > 0$ and solution
      of \eqref{eq:problem:sipoptbarrier}, and such that
      \begin{equation}
      \lim_{\mu\to 0} w(\mu, p) = w^\star\; .
      \end{equation}
    \item For $\mu$ small enough,
      \begin{equation}
        \label{eq:sensitivity:sipopt}
        J_p w(\mu, p) = - (J_w F_\mu)^{-1} J_p F_\mu \; .
      \end{equation}
  \end{itemize}
\end{proposition}
To solve the sparse linear system~\eqref{eq:sensitivity:sipopt}, sIpopt uses
the factorization of the KKT system computed by Ipopt at the local solution $w(\mu ,p)$.
To handle change in the active set, sIpopt uses a fix-and-relax strategy, without
resorting to the resolution of the QP~\eqref{eq:qpkyparisis}. The new constraints
activated are incorporated in the linear system~\eqref{eq:sensitivity:sipopt},
and the system is solved using a Schur-complement approach.

\paragraph{Variant 1: CasADi}
CasADi~\cite{andersson2019casadi} is a powerful modeler used to solve nonlinear programs~\eqref{eq:problem}.
The sensitivity analysis implemented in CasADi is slightly different than sIpopt~\cite{andersson2018sensitivity}.
Once a solution returned by the optimization solver, CasADi applies a primal-dual active
set method to determine exactly the active set at the solution (an
information that interior-point does not return by default, unless a crossover is used).
Then, CasADi uses a custom sparse QR solver to solve a non-symmetric formulation of
the KKT condition. Such QR factorization allows to quickly detect and resolve
the singularities associated to potential degeneracy.

\paragraph{Variant 2: QPTH}
sIpopt's method has been recently specialized in \cite{amos2017optnet} to the parallel extraction of the sensitivity
for a quadratic program (QP), for a fixed number of parameters $(p_1, \cdots, p_N)$.
The solver QPTH solves $N$ QP problems
in parallel using a specialized interior-point method leveraging batched dense linear algebra on the GPU.
Once a solution found, the sensitivity are evaluated using a formula analogous to \eqref{eq:sensitivity:sipopt}.
This has been proven relevant for machine learning applications, where neural network
are trained with batched stochastic gradient descent.
To the best of our knowledge, the method has not been extended to sparse problems, limiting
its adoption beyond machine learning. Also, little emphasis is put on the degenerate case.

\subsection{Path following method}
In nonlinear programming, one of the main application for sensitivity analysis is updating
the solution of \nlp\ after a small change in the parameter $p$:
once the directional derivative is obtained, we can compute an approximation of
a new primal solution $x(p')$ by setting the direction $h = p' - p$
and using the Taylor expansion~\cite{zavala2009advanced}:
\begin{equation}
  \label{eq:updateparam}
  x(p') = x(p) + x'(p, h) + o(\|h\|) \; .
\end{equation}
However, the method is valid only for small perturbations around~$p^\star$, without any active set change.
In case active set changes are occurring, they have to be accommodated locally using a path following method,

A path following method tracks the path of optimal solution between $p$ and the new parameter $p'$.
To do so, it defines the affine interpolation $p(t) = p + t(p' - p)$
for $t \in [0, 1]$ and tracks the solution of \nlp\ between $t=0$ and $t=1$.
This approach returns more accurate approximations, notably when facing active set changes.
For a sequence of scalars $t_{(0)} < t_{(1)} < \cdots < t_{(N)}$ such that $t_{(0)} = 0$
and $t_{(N)} = 1$, the perturbation between $t_{(m)}$ and $t_{(m+1)}$ is defined as
\begin{equation}
  h_{(m)} = (t_{(m+1)} - t_{(m)}) (p' - p) \;, \quad \forall m = 0, \cdots, N-1 \; .
\end{equation}
Then, we apply $N-1$ sensitivity updates with the direction $h_{(m)}$ to follow the path of the primal-dual solution $w(t)$
between $t_{(0)}$ and $t_{(N)}$.

We refer to \cite{zavala2009advanced} for a description of the path following
method. The method has been extended later on \cite{jaschke2014fast} for degenerate problems where
the dual solution is non-unique, using the method of Ralph and Dempe described in \S\ref{sec:solutionsensitivity:degenerate}.
The algorithm proposed in \cite{jaschke2014fast} tracks the active set changes
by updating iteratively the dual multipliers using the LP \eqref{eq:dempe:LP}.
As such, it can be interpreted as running $N$ iterations of a sequential quadratic programming (SQP) algorithm.
The original method~\cite{jaschke2014fast} didn't come with any convergence guarantee,
but it has been extended later in \cite{kungurtsev2017predictor} with a rigorous convergence
proof. In \cite{kungurtsev2017predictor}, the path following method is interpreted
as a predictor-corrector method, where the corrector solves a linear system to update
the sensitivity and the predictor solves a QP to track the active set changes. The multipliers
are again updated by solving a LP equivalent to \eqref{eq:dempe:LP}. We note
that the method requires the evaluation of the second-order derivatives of each constraint
in order to solve the QP in the predictor step.
To the best of our knowledge, we are not aware of a generic package that implements
the method described in \cite{kungurtsev2017predictor}.

\section*{Acknowledgements}
We thank Léonard Moracchini for his help in establishing the first
version of this manuscript, and Andrew Rosemberg and Tristan
Rigaut for their helpful feedbacks.

\appendix

\section{Derivatives}
Let $f: \Omega \subset E \to F$ a function, with $E$ and $F$ two normed spaces.
We start to recall the usual derivatives, before introducing the Clarke Jacobian more
adapted in the nonsmooth context.

\subsection{Differentiability}

\begin{definition}[Directional derivatives]
  Let $x \in \Omega$. The function $f$ has a \emph{directional derivative}
  at $x$ in the direction $h \in E$ if $x + th \in \Omega$ for $t$ small enough,
  and if the following limit exists
  \begin{equation}
    f'(x; h) = \lim_{t \downarrow 0} \; \dfrac{f(x + th) - f(x)}{t} \; .
  \end{equation}
\end{definition}

\begin{definition}[Gâteaux-differentiability]
  The function $f$ is Gâteaux differentiable (or $G$-differentiable) in $x \in \Omega$
  if there exists a directional derivatives for every direction $h \in E$,
  and if the application
  \begin{equation}
    h \to f'(x; h)  \; ,
  \end{equation}
  is \emph{continuous} and \emph{linear}.
  We note $f'(x)$ the operator such that $f'(x; h) = f'(x) \cdot h$.
\end{definition}
We note that a $G$-differentiable function is not necessarily continuous.

\begin{definition}[Fréchet-differentiability]
  The function $f$ is Fréchet differentiable ($F$-differentiable, or simply differentiable)
  at $x \in \Omega$ if there exists a bounded linear operator $L$ such that
  \begin{equation}
    \label{eq:deriv:frechet}
  \lim_{\|h\| \downarrow 0} \; \dfrac{1}{\|h\|}\Big(f(x + h) - f(x) - L h \Big) = 0 \; .
  \end{equation}
\end{definition}
$F$-differentiability implies that $f$ is continuous.
The operator $L$ is the \emph{derivative} of $f$ in $x$.
The condition \eqref{eq:deriv:frechet} rewrites equivalently
\begin{equation}
  f(x + h) = f(x) + L h + o(\|h\|) \; .
\end{equation}

The Hadamard derivative gives the directional derivative
along a curve tangential to $h$. It is required to study the sensitivity
of infinite-dimensional problems, as often encountered in statistics or
stochastic optimization.
\begin{definition}[Hadamard-differentiability]
  The function $f$ is Hadamard differentiable at $x$
  if for any mapping $\varphi: \mathbb{R}_+ \to X$ such that $\varphi(0) = x$
  and $\frac{1}{t} (\varphi(t) - \varphi(0))$ converges to a vector
  $h$ as $t \downarrow 0$, the limit
  \begin{equation}
    d_H f(x; h) = \lim_{t \downarrow 0}  \; \frac{1}{t} (f(\varphi(t)) - f(x)) \; ,
  \end{equation}
  does exist.
\end{definition}
For any sequences $\{h_n\}_n$ and $\{t_n\}_n$ such that $h_n \to h$
and $t_n \to 0^+$, the Hadamard directional derivative can also be written
in the form
\begin{equation}
  d_H  f(x; h) = \lim_{n\to\infty} \; \frac{1}{t_n} (f(x+t_n h_n) - f(x))  \; .
\end{equation}
The Hadamard derivative has the advantage of being compatible with the chain-rule
with minimal assumptions.
\begin{proposition}[Chain-rule]
  Let $f$ be Hadamard directionally differentiable at $x$,
  and $g$ be Hadamard directionally differentiable at $f(x)$. Then
  the composite mapping $g \circ f$ is Hadamard directionally differentiable
  at $x$, and satisfies the chain-rule $d_H (g \circ f)(x; h) = d_H g (f(x); d_H f(x; y))$.
\end{proposition}
Whereas the Hadamard derivative allows for different rates for every directions,
the Fréchet derivative imposes the rate is the same for each direction $h$.
For finite-dimensional space $E$, if the directional derivative is continuous
then the Fréchet and the Hadamard derivatives coincide.

\subsection{Generalized differentiability}
\label{sec:deriv:nonsmooth}
The notion of derivatives can be extended to nonsmooth functions.
We recall the following important results.

\begin{theorem}[Rademacher Theorem]
  Let $f: \Omega \to F$ be a locally Lipschitz continuous function.
  Then $f$ is differentiable almost everywhere in the sense of the
  Lebesgue measure.
\end{theorem}

The concept of Dini-derivatives generalizes the notion
of directional differentiability for non-smooth functions.
\begin{definition}[Dini-derivatives]
  Let $x \in \Omega$. The upper Dini derivative in the
  direction $h$ is defined as
  \begin{equation}
    D^+ f(x; h) = \lim_{t \downarrow 0} \sup \left(
      \dfrac{f(x + th) - f(x)}{t}
    \right) \;.
  \end{equation}
  Accordingly, the lower Dini derivative is defined as
  \begin{equation}
    D^- f(x; h) = \lim_{t \downarrow 0} \inf \left(
      \dfrac{f(x + th) - f(x)}{t}
    \right) \;.
  \end{equation}
\end{definition}

\begin{definition}[Bouligand-differential]
  Let $\mathcal{D}_f$ be the set of points at which $f$ is (Fréchet) differentiable.
  The $B$-differential of $f$ at $x \in \Omega$ is the set
  defined as
  \begin{equation}
    \partial_B f(x) = \big\{ J \in \mathcal{L}(E, F) \; : \;
      \exists \{x_k \}_k \in \mathcal{D}_f \; \text{such that} \; x_k \to x \;, f'(x_k) \to J
    \big\} \; .
  \end{equation}
\end{definition}

\begin{definition}[Clarke-differential]
  The $C$-differential of $f$ at $x \in \Omega$ is the convex
  hull of the Bouligand-differential, that is
  \begin{equation}
    \partial_C f(x) = \chull \; \partial_B f(x) \;.
  \end{equation}
\end{definition}

\begin{definition}[Clarke-generalized directional derivative] \\
  The Clarke-generalized directional derivative of $f$
  in the direction $h$ is defined by
  \begin{equation}
    d_C f(x; h) = \lim_{y \to x} \sup_{t \downarrow 0} \frac 1t \big[f(y + th) - f(y) \big] \; .
  \end{equation}
\end{definition}

\begin{proposition}
  Suppose that $f$ is $L$-Lipschitz near $x \in \Omega$. Then
  \begin{enumerate}
    \item $\partial_C f(x)$ is nonempty compact and convex,
    \item $\partial_C f(x)$ is locally bounded and upper semi-continuous at $x$.
    \item $d_C f(x; \cdot)$ is the support function of $\partial_C f(x)$, that is,
      \begin{equation}
        d_C f(x; h) = \max \{ s^\top h \; : \; \forall s \in \partial_C f(x) \}
        \;.
      \end{equation}
  \end{enumerate}
\end{proposition}

\begin{proposition}
  If $f$ is $G$-differentiable at $x$, then $f'(x) \in \partial_C f(x)$.
  If in addition $f$ is $C^1$ at $x$, then $\partial_C f(x) = \{ f'(x) \}$.
\end{proposition}

\begin{proposition}[Chain-rule]
  Suppose that $f$ is $L$-Lipschitz continuous on $\Omega$.
  If the function $g:F \to G$ is Lipschitz-continuous near $f(x)$,
  then
  \begin{equation}
    \partial_C (g \circ f) (x) \subset \chull \big\{ G F \; : \;
    G \in \partial_C g(f(x)) \;, \; F \in \partial_C f(x) \big\}
    \; .
  \end{equation}
\end{proposition}

\subsection{Semi-smooth functions}
\label{sec:deriv:semismooth}

\begin{definition}[Semi-smoothness]
  The function $f:\Omega \to F$ is semi-smooth in $x \in \Omega$ if:
  \begin{enumerate}
    \item $f$ is $L$-Lipschitz near $x$;
    \item $f$ has directional derivatives at $x$ in all directions;
    \item When $h \to 0$, $\sup_{J \in \partial_C f} \|f(x+h) -f(x) - Jh \| = o(\|h\|)$.
  \end{enumerate}
\end{definition}

\begin{proposition}
  Suppose $f$ is $L$-Lipschitz near $x$, and admits directional
  derivatives at $x$ in all directions. Then the following statements
  are equivalent:
  \begin{enumerate}
    \item $f$ is semi-smooth ;
    \item for $h \to 0$, $\sup_{J \in \partial_C f} \|Jh - f'(x;h)\| = o(\|h\|)$.
    \item for $h \to 0$ such that $x+h \in \mathcal{D}_f$, $f'(x+h) h - f'(x; h)= o(\|h\|)$.
  \end{enumerate}
\end{proposition}

\subsection{PC$^1$ functions}
\label{sec:deriv:pc1}
The class of $PC^1$ functions is a specific tool in nonsmooth analysis,
often employed for the sensitivity analysis of nonlinear programs.

\begin{definition}
  The function $f:X \to F$ is piecewise-differentiable (PC$^1$) at $x \in X$
  if there exists a neighborhood $N \subset X$ and a finite collection
  of $C^1$ selection functions $\mathcal{F}_f(x) = \{f_{(1)}, \cdots, f_{(k)} \}$ defined on $N$
  such that $f$ continuous on $N$ and for all $y \in N$, $f(y) \in
\{f_{(i)}(y)\}_{i \in [k]}$.
\end{definition}
If $f$ is PC$^1$, we define the active set associated at $x \in X$ as
\begin{equation}
  \activeset_f^{ess}(x) := \{i \in [k] \; : \; \forall y \in N , \; f(y) = f_{(i)}(y) \}\; .
\end{equation}

\begin{proposition}
  Let $f:X \to F$ a PC$^1$ function at $x \in X$. Then
  \begin{itemize}
    \item $\activeset_f^{ess}(x)$ is a non-empty set;
    \item $f$ is Lipschitz continuous near $x$ and directionally
      differentiable at $x$, with
      \begin{equation}
        \partial_B f(x) = \{ \nabla f_{(i)}(x) \; : \; i \in \activeset_f^{ess}(x) \}.
      \end{equation}
      In addition, the directional derivative $f'(x; \cdot)$ is piecewise linear.
  \end{itemize}
\end{proposition}

The notion of \emph{coherent orientation} plays a central role in obtaining
the invertibility of a $PC^1$ function in order to apply the implicit function
theorem.

\begin{definition}[Coherent orientation]
  Let $U \subset \mathbb{R}^n$ and $V \subset \mathbb{R}^m$ be open.
  A $PC^1$ function $F:U \times W \to \mathbb{R}^n$ is said
  \emph{coherently oriented} with respect to $x$ at $(x^\star, y^\star) \in U \times V$
  if all matrices in $\partial_{B,x} f(x^\star, y^\star)$ have the
  same non-vanishing determinant sign.
\end{definition}

\subsection{Lexicographic derivatives}
Lexicographic derivatives have been introduced by
Nesterov in \cite{nesterov2005lexicographic} as the derivatives
of \emph{lexicographically smooth} (L-smooth) functions.

\begin{definition}[Lexicographically smooth function]
  Let $f:X \to \mathbb{R}^m$ a locally Lipschitz continuous function on $X$.
  $f$ is \emph{lexicographically smooth} at $x$ if for any $k \in \mathbb{N}$
  and any $M = [m_{(1)}, \cdots, m_{(k)}] \in \mathbb{R}^{n \times k}$
  the following higher-order derivatives are well defined
  \begin{equation}
    \begin{aligned}
      & f^{(0)}_{x, M}: \mathbb{R}^n \to \mathbb{R}^m : d \to f'(x; d) \\
      & f^{(j)}_{x, M}: \mathbb{R}^n \to \mathbb{R}^m : d \to [f^{(j-1)}_{x, M}]'(m_{(j)}; d) & \forall j =1,\cdots, k \\
    \end{aligned}
  \end{equation}
\end{definition}
The class of L-smooth functions is closed under composition. $C^1$ functions,
convex functions and $PC^{1}$ functions are all L-smooth functions.

\begin{definition}[Lexicographic derivative]
  Let $f:X \to \mathbb{R}^m$ a L-smooth function at $x$, and $M \in \mathbb{R}^{n \times n}$ a nonsingular
  matrix. The \emph{lexicographic derivative} of $f$ at $x$ in the direction $M$
  is defined as
  \begin{equation}
    d_L f(x; M) = \nabla f^{(n)}_{x, M}(0)  \in \mathbb{R}^{m \times n} \;.
  \end{equation}
  The \emph{lexicographic subdifferential} of $f$ at $x$ is defined as
  \begin{equation}
    \partial_L f(x) = \{ d_L f(x; N) \; : \; N \in \mathbb{R}^{n\times n} \;,\; \text{det}(N) \neq 0 \}
    \;.
  \end{equation}
  For any $k \in \mathbb{N}$, $M = [m_{(1)}, \cdots, m_{(k)}] \in \mathbb{R}^{n \times k}$
  the \emph{LD-derivative} of $f$ in the directions $M$ is defined as
  \begin{equation}
    f'(x; M) = \big[
      f^{(0)}_{x, M}(m_{(1)}), ~
      f^{(1)}_{x, M}(m_{(2)}), ~
      \cdots,
      f^{(k-1)}_{x, M}(m_{(k)})
    \big]  \; .
  \end{equation}
\end{definition}
The LD-derivative is uniquely defined for $M \in \mathbb{R}^{n \times k}$.
If in addition $M \in \mathbb{R}^{n \times n}$ is square and
nonsingular, it satisfies
\begin{equation}
  f'(x; M) = d_L f(x; M) M \; .
\end{equation}
Interestingly, the LD-derivative obeys a sharp chain-rule.
\begin{proposition}[Chain-rule]
  Let $h:X \to Y$ and $g:Y \to \mathbb{R}^q$ be L-smooth at $x \in X$
  and $h(x) \in Y$, respectively. Then the function $g\circ h$
  is L-smooth at $x \in X$, and for any $k \in \mathbb{N}$
  and $M \in \mathbb{R}^{n \times k}$,
  \begin{equation}
    [g \circ h]'(x; M) = g'(h(x); h'(x; M)) \; .
  \end{equation}
\end{proposition}

LD-derivatives give convenient calculus rules for the derivatives
of $\min$, $\max$, abs, and other nonsmooth functions widely encountered in practice.
These rules are based on \emph{lexicographic ordering}, explaining why LD-derivatives
are called "lexicographic".
For given vectors $x, y \in \mathbb{R}^n$, we say that $x$ is lexicographically
lower than $y$ if
\begin{equation}
  \begin{aligned}
  & x \prec y ~ \text{if and only if} ~ \exists j \in [n] ~ \text{such that} ~ x_i = y_i ~ \forall i < j ~
  \text{and} ~ x_j < y_j \; , \\
  & x \preceq y ~ \text{if and only if} ~ x = y ~ \text{or} ~ x \prec y \; .
  \end{aligned}
\end{equation}
The \emph{lexicographic minimum function} returns the lexicographically ordered
minimum of two vectors:
\begin{equation}
  \textbf{Lmin}: \mathbb{R}^n \times \mathbb{R}^n \to \mathbb{R}^n:
  (x, y) \mapsto \left\{\begin{aligned}
      & x & \text{if} ~ x \preceq y \; ,  \\
      & y & \text{if} ~ x \succ y \; .
  \end{aligned}
  \right.
\end{equation}
Similarly, the \emph{lexicographic matrix minimum function} compares two matrices
row by row, and is defined as
\begin{equation}
  \label{eq:deriv:lmmin}
  \textbf{LMmin}: \mathbb{R}^{m \times n} \times \mathbb{R}^{m \times n} \to \mathbb{R}^{m \times n}:
  (X, Y) \mapsto \begin{bmatrix}
    \textbf{Lmin}(X_1^\top, Y_1^\top) \\
    \textbf{Lmin}(X_2^\top, Y_2^\top) \\
    \cdots \\
    \textbf{Lmin}(X_m^\top, Y_m^\top) \\
  \end{bmatrix} \; .
\end{equation}

\section{Multivalued mapping}
A multivalued mapping (or point-to-set map) is a mapping
$\Gamma: E \rightrightarrows F$ where for each $x \in E$,
$\Gamma(x)$ is a subset of $F$.
They have been used extensively in mathematical programming
to study the solution maps of optimization problems.

\begin{definition}
  Let $\Gamma: E \rightrightarrows F$ a point-to-set map.
  \begin{itemize}
    \item $\Gamma$ is \emph{upper semicontinuous}
      at $x \in E$ if for each open set $\Omega \subset F$ satisfying $\Gamma(x) \subset \Omega$
      there exists a neighborhood $N(x)$ of $x$ such that for all $y \in N(x)$,
      $\Gamma(y) \subset \Omega$.
    \item $\Gamma$ is \emph{lower semicontinuous}
      at $x \in E$ if for each open set $\Omega \subset F$ satisfying $\Gamma(x) \cap \Omega \neq \varnothing$
      there exists a neighborhood $N(x)$ of $x$ such that for all $y \in N(x)$,
      $\Gamma(y) \cap \Omega \neq \varnothing$.
    \item $\Gamma$ is \emph{continuous} at $x$ if it is lower semicontinuous and
        upper semicontinuous at $x$.
  \end{itemize}
\end{definition}

\begin{definition}
  The point-to-set map $\Gamma:E \rightrightarrows F$ is \emph{closed}
  at $x$ if there exists a sequence $x_n \in E$, $x_n\to x$ such that
  $y_n \in \Gamma(x_n)$ and $y_n \to y$ imply $y \in \Gamma(x)$.
\end{definition}
\begin{definition}
  The point-to-set map $\Gamma:E \rightrightarrows F$ is \emph{open}
  at $x$ if $y \in \Gamma(x)$ implies there exists a sequence $x_n \in E$, $x_n\to x$ such that
  there exists $m$ and $\{y_n\}$ such that $y_n \in \Gamma(x_n)$ for all $n \geq m$
  and $y_n \to y$.
\end{definition}

\begin{definition}
  The point-to-set map $\Gamma:E \rightrightarrows F$ is \emph{uniformly compact}
  near $\hat{x}$ if the set $\cup_{x \in N(\hat{x})} \Gamma(x)$ is bounded
  for some neighborhood $N(\hat{x})$ of $\hat{x}$.
\end{definition}
If $\Gamma$ is uniformly compact near $\hat{x}$, then $\Gamma$ is closed
if and onl if $\Gamma(\hat{x})$ is a compact set and $\Gamma$ is upper semicontinuous
at $\hat{x}$.

\begin{definition}
  The point-to-set map $\Gamma: E \rightrightarrows F$ is \emph{upper Lipschitzian}
  with modulus $L$ at $\hat{x} \in E$ if there is a neighborhood
  $N$ of $\hat{x}$ such that for each $x \in N$,
  $\Gamma(x) \subset \Gamma(\hat{x}) + L \|x - \hat{x} \|B$,
  with $B$ the unit-ball in $F$.
\end{definition}

\section{A primer on generalized equations}
\label{sec:robinson}
Let $X$ be a linear space,
and $C$ be a non-empty closed convex set in $X$.

\begin{definition}[Normal cone]
  We define the \emph{normal cone} operator associated to the set $C$
  as, for $x \in X$
  \begin{equation}
    N_C(x) =
    \left\{
    \begin{aligned}
    & \{ y \in X^\star \; : \; \langle y , v - x \rangle \leq 0  \;, \; \forall v \in C \} & \text{if} ~ x \in C \; ,\\
    & \varnothing\;, & \text{if}~ x \notin C \; .
    \end{aligned}
    \right.
  \end{equation}
\end{definition}
With the normal cone, we can define formally the notion of \emph{variational inequality},
generalizing optimization problem.

\begin{definition}[Variational inequality]
  Let $p \in P$ be a parameter. For $f: X \times P \to X^\star$, we define
  the \emph{paramaterized variational inequality} (also called \emph{generalized
  equation}) the problem of finding $x \in X$ such that
  \begin{equation}
    \label{eq:varineq}
    f(x, p) + N_C(x) \ni 0 \; .
  \end{equation}
\end{definition}
The equation \eqref{eq:varineq} is a compact formulation stating
that $-f(x, p) \in N_C(x)$, or, equivalently
\begin{equation}
  \langle f(x, p) , v - x \rangle \geq 0 \; , \quad \forall v \in C \; .
\end{equation}
The solution mapping associated to \eqref{eq:varineq} is defined as
\begin{equation}
  S(p) = \{ x \in X \; : \; f(x, p) + N_C(x) \ni 0 \}  \; .
\end{equation}
Robinson extended the Implicit function theorem to the nonsmooth
variational inequality \eqref{eq:varineq}, by extending the notion
of \emph{regular solution} to the variational inequality setting,
where the solution mapping $S(p)$ is not necessarily single-valued.
Robinson introduced the notion of \emph{strong regularity}, which
extends the non-singularity condition we impose in the implicit
function theorem. If we suppose that $f$ is (Fréchet) differentiable,
the idea is to look at the linearized generalized
equation around a solution $x_0$
\begin{equation}
  T(x, p) \ni 0
  \; ,
\end{equation}
where we have defined the linearized operator at $x_0$:
\begin{equation}
T(x, p) := f(x_0, p) + \nabla_x f(x_0, p) (x - x_0) + N_C(x) \; .
\end{equation}

\begin{definition}[Strong regularity~\cite{robinson1980strongly}]
  \label{def:strongregular}
  We say that \eqref{eq:varineq} is \emph{strongly regular} at
  a solution $x_0$ with Lipschitz constant $L$ if
  there exists a neighborhood $U \subset X^\star$ of
  the origin and $V \subset X$ of $x_0$ such that the restriction
  to $U$ of $T^{-1} \cap V$ is a \emph{single-valued} function
  from $U$ to $V$ and is $L$-Lipschitzian on $U$.
\end{definition}
Putting it more simply, strong regularity implies that the
inverse of the linearized operator $T^{-1}$ has a Lipschitz continuous
\emph{single-valued} localization at $0$ for $x_0$.
If $C$ is the whole-space $X$, then strong regularity implies that
the Jacobian $\nabla_x f(x_0, p)$ is non-singular.
The following theorem extends the implicit function theorem to
generalized equations~\eqref{eq:varineq} satisfying the strong regularity condition
\ref{def:strongregular}.
\begin{theorem}[Theorem 2.1, \cite{robinson1980strongly}]
  \label{thm:robinsonift}
  Let $p \in P$ be a given parameter, $x_0$ solution of the
  generalized equation~\eqref{eq:varineq}.
  Suppose that $f$ is Fréchet-differentiable w.r.t. $x$ exists at $p$,
  and that both $f(\cdot, \cdot)$ and $\nabla_x f(\cdot, \cdot)$ are
  continuous at $(x_0, p)$. If \eqref{eq:varineq} is \emph{strongly regular}
  at $x_0$, then for any $\varepsilon >0$ there exists neighborhood $N$ of $p$ and $W$ of $x_0$
  as well as a single-valued function $x:N \to W$ sucht that for any $p \in N$,
  $x(p)$ is the unique solution of the inclusion
  \begin{equation}
    f(x, p) + N_C(x) \ni 0 \; ,
  \end{equation}
  and, for all $p, q \in N$, one has
  \begin{equation}
    \| x(p) - x(q) \| \leq (L + \varepsilon) \| f(x(p), p) - f(x(q), q) \|
    \; .
  \end{equation}
\end{theorem}
Theorem~\ref{thm:robinsonift} states that the strong regularity
of \eqref{eq:varineq} implies that the solution mapping $S(p)$
of \eqref{eq:varineq} has a \emph{Lipschitz-continuous} single-valued
localization at $p$ for $x_0$.

\begin{corollary}[Locally-Lipschitz condition]
  Suppose the hypotheses of Theorem~\ref{thm:robinsonift} hold. Suppose in
  addition there exists a constant $\kappa$ such that for each $(p, q) \in N$,
  for each $x \in W$, one has
  \begin{equation}
    \| f(x, p) - f(x, q) \| \leq \kappa \| p - q \| \; .
  \end{equation}
  Then $x(\cdot)$ is Lipschitzian on $N$ with modulus $\kappa (L + \varepsilon)$.
\end{corollary}
\end{document}